\documentclass[12pt,reqno]{amsart}
\usepackage{enumerate, latexsym, amsmath, amsfonts, amssymb, amsthm, color}
\def\pmod #1{\ ({\rm{mod}}\ #1)}
\def\Z{\Bbb Z}
\def\N{\Bbb N}

\def\l{\left}
\def\r{\right}
\def\bg{\bigg}
\def\({\bg(}
\def\){\bg)}
\def\t{\text}
\def\f{\frac}

\def\ls{\leqslant}
\def\gs{\geqslant}
\def\se {\subseteq}
\def\sm{\setminus}

\def\ve{\varepsilon}

\def\eq{\equiv}

\def\Ack{\medskip\noindent {\bf Acknowledgments}}
\theoremstyle{plain}

\newtheorem{conjecture}{Conjecture}
\theoremstyle{definition}

\theoremstyle{remark}
\newtheorem{remark}{Remark}

 \vspace{4mm}

\begin{document}

\hbox{In: M. Nathanson (ed.), {\it Combinatorial and Additive Number Theory II:}}
\hbox{CANT, New York, NY, USA, 2015 and 2016, Springer Proc. in Math. \& Stat.,}
\hbox{Vol. 220, Springer, New York, 2017, pp. 279-310.}
\medskip

\title
[{Conjectures on representations involving primes}]
{Conjectures on representations\\ involving primes}

\author
[Zhi-Wei Sun] {Zhi-Wei Sun}

\address {Department of Mathematics, Nanjing
University, Nanjing 210093, People's Republic of China}
\email{zwsun@nju.edu.cn}

\keywords{Conjectures, primes, practical numbers, representations.
\newline \indent 2010 {\it Mathematics Subject Classification}. Primary  11A41, 11P32; Secondary 11B13, 11D68.
\newline \indent The initial version of this paper was posted to {\tt arXiv} in Nov. 2012 as a preprint with the ID {\tt arXiv:1211.1588}.}

 \begin{abstract} We pose 100 new conjectures on representations involving primes or related things, which might interest number theorists and stimulate further research.
Below are five typical examples:
(i) For any positive integer $n$, there exists $k\in\{0,\ldots,n\}$ such that $n+k$ and $n+k^2$ are both prime.
(ii) Each integer $n>1$ can be written as $x+y$ with $x,y\in\{1,2,3,\ldots\}$ such that $x+ny$ and $x^2+ny^2$ are both prime.
(iii) For any rational number $r>0$, there are distinct primes $q_1,\ldots,q_k$ with $r=\sum_{j=1}^k1/(q_j-1)$.
(iv) Every $n=4,5,\ldots$ can be written as $p+q$, where $p$ is a prime with $p-1$ and $p+1$ both practical, and $q$ is either prime or
practical.
(v) Any positive rational number can be written as $m/n$, where $m$ and $n$ are positive integers with $p_m+p_n$ a square (or $\pi(m)\pi(n)$
a positive square), $p_k$ is the $k$-th prime and $\pi(x)$ is the prime-counting function.
\end{abstract}

\maketitle

\section{Introduction}
\setcounter{lemma}{0}
\setcounter{theorem}{0}
\setcounter{corollary}{0}
\setcounter{remark}{0}
\setcounter{equation}{0}
\setcounter{conjecture}{0}

Primes have been investigated for over two thousand years. Nevertheless, there are many problems on primes remain open.
The famous Goldbach conjecture (cf. \cite{CP} and \cite{N}) states that any even integer $n>2$ can be represented as a sum of two primes.
Lemoine's conjecture (see \cite{L}) asserts that any odd integer $n>6$ can be written as $p+2q$ with $p$ and $q$ both prime;
this is a refinement of the weak Goldbach conjecture (involving sums of three primes) proved by I. M. Vinogradov \cite{V} for large odd numbers
and confirmed by H. A. Helfgott \cite{He} completely.
Legendre's conjecture states that for any
positive integer $n$ there is a prime between $n^2$ and $(n+1)^2$. Another well known conjecture of A. de Polignac asserts that
for any positive even number $d$ there are infinitely many positive integers $n$
with $p_{n+1}-p_n=d$, where $p_k$ denotes the $k$-th prime. (This conjecture in the case $d=2$
is the famous twin prime conjecture; recently Y. Zhang \cite{Z} made an important breakthrough along this line.)
Polignac's conjecture follows from the following well-known hypothesis due to A. Schinzel.
\medskip

\noindent {\bf Schinzel's Hypothesis}. {\it If
$f_1(x),\ldots,f_k(x)$ are irreducible polynomials with integer
coefficients and positive leading coefficients such that there is no
prime dividing the product $f_1(q)f_2(q)...f_k(q)$ for all $q\in\Z$,
then there are infinitely many positive integers $n$
 such that $f_1(n),f_2(n),\ldots,f_k(n)$ are all primes.}

 \medskip

A positive integer $n$ is said to be {\it practical} if every $m=1,\ldots,n$ can be written as the sum of some distinct (positive) divisors of $n$.
In 1954 B. M. Stewart \cite{St} showed that if $q_1<\cdots<q_r$ are distinct primes and $a_1,\ldots,a_r$ are positive integers
then $m=q_1^{a_1}\cdots q_r^{a_r}$ is practical if and only if $q_1=2$ and
$$q_{s+1}-1\ls\sigma(q_1^{a_1}\cdots q_s^{a_s})\quad\t{for all}\ \ 0<s<r,$$
where $\sigma(n)$ stands for the sum of all divisors of $n$.
The behavior of practical numbers is quite similar to that of primes.
For example, G. Melfi \cite{Me} proved the following Goldbach-type conjecture of M. Margenstern \cite{Ma}:
Each positive even integer is a sum of two practical numbers, and there are infinitely many practical numbers $m$ with $m-2$ and $m+2$ also practical.
Recently, A. Weingartner \cite{W} proved that the number of practical numbers not exceeding $x\gs2$ is asymptotically equivalent to $cx/\log x$, where
$c$ is a positive constant close to $1$; this analog of the Prime Number Theorem for practical numbers was first conjectured by Margenstern \cite{Ma} in 1991.

In the published papers \cite{S13a, S13b, S15, S17} the author posed many conjectures on primes. For example, \cite{S15} contains 60 problems on combinatorial properties
of primes many of which depend on some exact values of the prime-counting function $\pi(x)$. ($\pi(x)$ with $x\gs0$
denotes the number of primes not exceeding $x$.)

In this paper we present 100 new conjectures on representations involving primes or related things.
In particular, we find some surprising
refinements of Goldbach's conjecture, Lemoine's conjecture, Legendre's conjecture and the twin prime conjecture.
The next section contains 25 conjectures, the first of which is a general hypothesis (similar to Schinzel's Hypothesis) on representations of integers involving primes,
and the other 24 conjectures are closely related to this general hypothesis.
In Section 3 we include 45 conjectures on various other representation problems for integers.
In Section 4 we pose 30 conjectures on representations of positive rational numbers and related things.
For numbers of representations related to some conjectures in Sections 2-4, the reader may consult \cite{S} for certain sequences in the OEIS.

We hope that the 100 conjectures collected here might interest some number theorists and stimulate further research.

Throughout this paper, we set $\N=\{0,1,2,\ldots\}$ and $\Z^+=\{1,2,3,\ldots\}$.
For a real number $x$, the fractional part of $\{x\}$ is given by $x-\lfloor x\rfloor$.
For $a\in\Z$ and $n\in\Z^+$, by $\{a\}_n$ we mean the least nonnegative residue of $a$ modulo $n$, i.e., $\{a\}_n=n\{a/n\}$.
For $a\in\Z$ and $n\in\Z^+$ with $2\nmid n$, $(\f an)$ denotes the Jacobi symbol. As usual, $\varphi$ stands for Euler's totient function.

\section{A general hypothesis and related conjectures}
\setcounter{lemma}{0}
\setcounter{theorem}{0}
\setcounter{corollary}{0}
\setcounter{remark}{0}
\setcounter{equation}{0}
\setcounter{conjecture}{0}

Note that Schinzel's Hypothesis does not imply Goldbach's conjecture.
Here we pose a general hypothesis on representations of integers.

\begin{conjecture}\label{Conj2.1} {\rm (General Hypothesis, 2012-12-28)} Let
$$f_1(x,y),\ \ldots,\ f_m(x,y)$$ be non-constant polynomials with integer coefficients.
Suppose that for all large $n\in\Z^+$,  those $f_1(x,n-x),\ldots,f_m(x,n-x)$ are irreducible, and there is no prime dividing
all the products $\prod_{k=1}^mf_k(x,n-x)$ with $x\in\Z$. If $n\in\Z^+$ is large enough, then we can write $n=x+y\ (x,y\in\Z^+)$
such that $|f_1(x,y)|,\ldots,|f_m(x,y)|$ are all prime.
\end{conjecture}
\begin{remark}\label{Rem2.1} In view of this general hypothesis, almost all of the other conjectures in this section are essentially reasonable.
\end{remark}

\begin{conjecture}\label{Conj2.2} {\rm (Symmetric Conjecture, 2015-08-27)} For any integer $n>6$, there is a prime $p<n/n'$ such that $n-(pn'-1)$ and $n+(pn'-1)$ are both prime,
where $n'=2-\{n\}_2$ is $1$ or $2$ according as $n$ is odd or even.
\end{conjecture}
\begin{remark}\label{Rem2.2} Conjecture \ref{Conj2.2} is stronger than Goldbach's conjecture and Lemoine's conjecture.
We have verified Conjecture \ref{Conj2.2} for all $n=7,\ldots,10^8$, see \cite[A261627 and A261628]{S} for related data.
Conjecture \ref{Conj2.1} implies that Conjecture \ref{Conj2.2} holds for all sufficiently large integers $n$.
In fact, if we apply Conjecture \ref{Conj2.1} with $f_1(x,y)=x$, $f_2(x,y)=2y+1$ and $f_3(x,y)=4x+2y-1$, then for sufficiently large $n\in\Z^+$
there are primes $p$ and $q$ with $n=p+(q-1)/2$ (i.e., $2n-(2p-1)=q$) such that $2n+2p-1=4p+q-2$ is prime; if we apply Conjecture \ref{Conj2.1} with $f_1(x,y)=2x+1$, $f_2(x,y)=2y-1$ and $f_3(x,y)=4x+2y-1$, then for sufficiently large $n\in\Z^+$
there are primes $p$ and $q$ with $n=(p-1)/2+(q+1)/2$ (i.e., $2n-1-(p-1)=q$) such that $2n-1+(p-1)=2p+q-2$ is prime.
\end{remark}

\begin{conjecture}\label{Conj2.3} For each $n=6,7,\ldots$  there is a prime $p<n$ such that
both $6n-p$ and $6n+p$ are prime.
\end{conjecture}
\begin{remark}\label{Rem2.3} We also have some conjectures involving practical numbers similar to Conjectures \ref{Conj2.2} and \ref{Conj2.3}, see \cite[A261641]{S}
and Conjectures \ref{Conj3.43} and \ref{Conj3.44}. Conjecture \ref{Conj2.1} with $f_1(x,y)=x$, $f_2(x,y)=5x+6y$ and $f_3(x,y)=7x+6y$
implies that Conjecture \ref{Conj2.3} holds for sufficiently large integers $n$.
\end{remark}

\begin{conjecture}\label{Conj2.4} {\rm (2012-12-22)} Any integer $n\gs12$ can be written as $p+q\ (q\in\Z^+)$ with $p,\,p+6,\,6q-1$ and $6q+1$ all prime.
\end{conjecture}
\begin{remark}\label{Rem2.4} Conjecture \ref{Conj2.1} implies that Conjecture \ref{Conj2.4} holds for all sufficiently large integers $n$.
We have verified Conjecture \ref{Conj2.4} for $n$ up to $10^9$, see \cite[A199920]{S} for numbers of such representations.
Conjecture \ref{Conj2.4} implies that there are infinitely many twin primes and also infinitely many sexy primes, because
for any $m=2,3,\ldots$ the interval $[m!+2,m!+m]$ of length $m-2$ contains no prime.
\end{remark}

\begin{conjecture}\label{Conj2.5} {\rm (2013-10-09)} Any integer $n>1$ can be written as $k+m\ (k,m\in\Z^+)$ with $6k-1$ a Sophie Germain prime and $\{6m-1,6m+1\}$ a twin prime pair.
\end{conjecture}
\begin{remark}\label{Rem2.5} Recall that a Sophie Germain prime is a prime $p$ with $2p+1$ also prime.
Conjecture \ref{Conj2.1} implies that Conjecture \ref{Conj2.5} holds for all sufficiently large integers $n$.
We have verified  Conjecture \ref{Conj2.5} for all $n=2,\ldots,10^8$, see \cite[A227923]{S} for numbers of such representations.
Conjecture \ref{Conj2.5} implies that there are infinitely many twin primes and also infinitely many Sophie Germain primes.
For example, if all twin primes do not exceed an integer $N>2$ and $(N+1)!/6=k+m\ (k,m\in\Z^+)$ with $6k-1$ a Sophie Germain prime
 and $\{6m-1,6m+1\}$ a twin prime pair, then $6k-1=(N+1)!-(6m+1)$ with $2\ls 6m+1\ls N$ which contradicts that $6k-1$ is prime.
\end{remark}

Recall that for two subsets $X$ and $Y$ of $\Z$ the sumset $X+Y$ is defined as $\{x+y:\ x\in X\ \mbox{and}\ y\in Y\}$.

\begin{conjecture}\label{Conj2.6} {\rm (2013-01-03)} Let
\begin{align*} A=&\{x\in\Z^+:\
6x-1\ \t{and}\ 6x+1\ \t{are both prime}\},
\\ B=&\{x\in\Z^+:\ 6x+1\ \t{and}\ 6x+5\ \t{are both prime}\},
\\ C=&\{x\in\Z^+:\ 2x-3\ \t{and}\ 2x+3\ \t{are both prime}\}.
\end{align*}
Then
$$A+B=\{2,3,\ldots\},\ B+C=\{5,6,\ldots\},\ \ A+C=\{5,6,\ldots\}\sm\{161\}.$$
Also, if we set $2X:=X+X$ for $X\se\Z$, then
$$2A\supseteq\{702,703,\ldots\},\ 2B\supseteq\{492,493,\ldots\},\
2C\supseteq\{4006,4007,\ldots\}.$$
\end{conjecture}
\begin{remark}\label{Rem2.6} Conjecture \ref{Conj2.1} implies that each of the sumsets $A+B, B+C, A+C, 2A, 2B, 2C$ in Conjecture \ref{Conj2.6} contain all sufficiently large integers.
\end{remark}

\begin{conjecture}\label{Conj2.7} {\rm (2013-10-12)} {\rm (i)} For any integer $n>3$, we can write $2n$ as $p+q$ with $p,q,3p-10,3q+10$ all prime.

{\rm (ii)} For any integer $n>4$ not equal
to $76$, we can write $2n$ as $p+q$ with $p,3p-10,q,3q-10$ all prime.
\end{conjecture}
\begin{remark}\label{Rem2.7} Note that if $2n=p+q$ then $6n=(3p-10)+(3q+10)$. We have verified Conjecture \ref{Conj2.7} for $n$ up to $10^8$. See \cite[A230230]{S}
for related data. Conjecture \ref{Conj2.1} implies that Conjecture \ref{Conj2.7} holds for all sufficiently large integers $n$.
\end{remark}

\begin{conjecture}\label{Conj2.8} {\rm (2012-11-07)} For any integer $n>8$, we can write $2n-1$ as $p+2q$ with $p,q$ and $p^2+60q^2$ all prime.
\end{conjecture}
\begin{remark}\label{Rem2.8}  This is stronger than Lemoine's conjecture. We have verified Conjecture \ref{Conj2.8} for $n$ up to $10^8$. See
\cite[A218825]{S} for related data.
\end{remark}

\begin{conjecture}\label{Conj2.9} {\rm (2013-10-16)}  Any integer $n>3$ can be written
as $p+q\ (q\in\Z^+)$ with $p$, $2p^2-1$ and $2q^2-1$ all prime.
\end{conjecture}
\begin{remark}\label{Rem2.9} See \cite[A230351]{S} for related data.
Note that each of 7, 12, 68, 330 has a unique required
representation:
\begin{gather*} 7 = 3+4,\  2\cdot3^2-1=17,\ 2\cdot4^2-1=31;
\\12=2+10,\ 2\cdot2^2-1=7,\ 2\cdot10^2-1=199;
\\68=43+25,\ 2\cdot43^2-1=3697,\ 2\cdot25^2-1=1249;
\\330= 7+323,\ 2\cdot7^2-1=97,\ 2\cdot323^2-1=208657.
\end{gather*}
\end{remark}

In 2001 A. Murthy (cf. \cite{Mu}) conjectured that for any integer $n>1$ there is
an integer $0<k<n$ such that $kn+1$ is prime. In 2005 he \cite{Mu} conjectured
any integer $n>3$ can be written as $x+y\ (x,y\in\Z^+)$ with $xy-1$ prime.
In 1990s Ming-Zhi Zhang (cf. \cite[p.\,161]{G}) asked whether any odd integer $n>1$
can be written as $a+b$ with $a,b\in\Z^+$ and $a^2+b^2$ prime.

\begin{conjecture}\label{Conj2.10} {\rm (2012-12-20) (i)} For any integer $n>3$, there is an integer $k\in\{1,\ldots,n-1\}$
such that $kn+1$ and $k(n-k)-1$ are both prime.

{\rm (ii)} For any odd integer $n>1$, there is an integer $k\in\{1,\ldots,n-1\}$
such that $kn+1$ and $k^2+(n-k)^2$ are both prime.
\end{conjecture}
\begin{remark}\label{Rem2.10} This combines Murthy's conjectures and Zhang's conjecture.
We also conjecture that any integer $n>3$ can be written as $x+y$ with $x,y\in\Z^+$ such that
$3x\pm1$ and $xy-1$ are all prime (cf. \cite[A220431]{S}).
\end{remark}

\begin{conjecture}\label{Conj2.11} {\rm (2013-11-12)} {\rm (i)} Any integer $n>2$ can be written
as $k+m\ (k,m\in\Z^+)$ with $k^2m-1$ prime. Also, each integer $n>4$ can be written
as $k+m\ (k,m\in\Z^+)$ with $k^2m+1$ prime.

{\rm (ii)} Any integer $n>1$ can be written
as $k+m\ (k,m\in\Z^+)$ with $(km)^2+km+1$ prime. Also, each integer $n>2$ can be written
as $k+m\ (k,m\in\Z^+)$ with $(km)^2+km-1\ (or\ 2k^2m^2-1)$ prime.
\end{conjecture}
\begin{remark}\label{Rem2.11} See \cite[A231633]{S} for related data.
\end{remark}

\begin{conjecture}\label{Conj2.12} {\rm (2013-10-13) (i)}
For any integer $n>1$, there
is a prime $p\ls n$ such that $(p-1)n+1$ is prime.
Moreover, for any integer $n>4$, there
is a prime $p<n$ such that $3p+8$ and $(p-1)n+1$ are both prime.

{\rm (ii)} Any integer $n>5$ can be written
as $p+q\ (q\in\Z^+)$ with $p,\, 3p-10$ and $(p-1)q-1$ all prime.
\end{conjecture}
\begin{remark}\label{Rem2.12} See \cite[A230243 and A230241]{S} for related data.
\end{remark}

\begin{conjecture}\label{Conj2.13} {\rm (2012-12-16)} For any integer $n>1$, we can write $2n$
as $p+q$, where $p$ is a Sophie Germain prime, $q$ is a positive integer, and $(p-1)^2+q^2$
is prime.
\end{conjecture}
\begin{remark}\label{Rem2.13} This is stronger than Zhang's conjecture. Conjecture \ref{Conj2.1} implies that any sufficiently large $n$
can be written as $x+y\ (x,y\in\Z^+)$ with $p=2x+1$, $2p+1=4x+3$ and $$(p-1)^2+(2n-p)^2=(2x)^2+(2y-1)^2$$ all prime.
See \cite[A220554]{S} for related data. For example, $32=11+21$ with $11$ a Sophie Germain prime and $(11-1)^2+21^2=541$ a prime.
\end{remark}

\begin{conjecture}\label{Conj2.14} {\rm (i) (2011-11-04)} Any odd integer $n>1$ can be written as $x+y$ with $x,y\in\Z^+$ such that $x^4+y^4$ is prime.

{\rm (ii) (2012-12-01)} Any integer $n>10$ can be written as $p+q\ (q\in\Z^+)$ with $p$, $p+6$ and $p^2+3pq+q^2=n^2+pq$ all prime.

{\rm  (iii) (2013-11-21)} Let $n>1$ be an odd integer. We can write $n=k+m$ with $k,m\in\Z^+$
such that both $k^2+m^2$ and $k^3+m^2$ are prime.
\end{conjecture}
\begin{remark}\label{Rem2.14} See \cite[A218656, A218654, A218754 and A232269]{S} for related data.
\end{remark}

\begin{conjecture}\label{Conj2.15} {\rm  (Olivier Gerard and Zhi-Wei Sun, 2013-10-13)}. For any integer $n>1$, we can write $2n$ as $p+q$ with $p,\,q$ and
$(p-1)(q+1)-1$ all prime.
\end{conjecture}
\begin{remark}\label{Rem2.15} This is stronger than Goldbach's conjecture.
Note also that $(p-1)+(q+1)=p+q$. We have verified Conjecture \ref{Conj2.15} for all $n=2,\ldots,10^8$. See \cite[A227909]{S} for related data.
 \end{remark}

\begin{conjecture}\label{Conj2.16} {\rm (2012-11-30)} Any integer $n>7$ can be
written as $p+q$ $(q\in\Z^+)$ with $p$ and $2pq+1$ both prime.
In general, for each $m\in\N$ any sufficiently large integer $n$ can be written as $x+y\ (x,y\in\Z^+)$ with $x-m,\,x+m$ and $2xy+1$ all prime.
\end{conjecture}
\begin{remark}\label{Rem2.16} We have verified the first assertion in Conjecture \ref{Conj2.16} for all $n=8,9,\ldots,10^9$. See \cite[A219864]{S} for related data.
Concerning the general statement in Conjecture \ref{Conj2.16}, for $m=1,2,3,4,5,6,7,8,9,10$ it suffices to require that $n$ is greater than
$$623,\, 28,\, 151,\,357,\,199,\,307,\,357,\,278,\, 697,\,263$$
respectively.
\end{remark}
\medskip

\begin{conjecture}\label{Conj2.17} {\rm (2013-10-14)} Any integer $n>3$ can be written
as $p+q\ (q\in\Z^+)$ with $p$ and $(p+1)q/2+1$ both prime.
\end{conjecture}
\begin{remark}\label{Rem2.17} We have verified this conjecture for all $n=4,\ldots,10^8$.
See \cite[A230254]{S} for related data. For example, 30 has a unique representation $2+28$ with
$(2+1)28/2+1=43$ prime.
\end{remark}

Bertrand's Postulate proved by Chebyshev in 1852 states that for any
positive integer $n$, the interval $[n,2n]$ contains at least a
prime. Goldbach's conjecture essentially asserts that for any integer $n>1$ there is
an integer $k\in\{0,\ldots,n\}$ such that $n-k$ and $n+k$ are both
prime. The following conjecture is of a similar flavor.

\begin{conjecture}\label{Conj2.18} {\rm (2012-12-18)} For each positive integer $n$,
there is an integer $k\in\{0,\ldots,n\}$ such that $n+k$ and $n+k^2$
are both prime.
\end{conjecture}
\begin{remark}\label{Rem2.18} We have verified this for $n$ up to $10^8$.
See \cite[A185636 and A204065]{S} for related data. The author would like to offer 100 US dollars as the prize for the first solution of Conjecture \ref{Conj2.18}.
\end{remark}

\begin{conjecture}\label{Conj2.19} {\rm (2013-04-15)} For any positive integer $n$,
there is a positive integer $k\ls 4\sqrt{n+1}$ such that $n^2+k^2$
is prime.
\end{conjecture}
\begin{remark}\label{Rem2.19} Note that the least $k\in\Z^+$ with $63^2+k^2$ prime
is $32=4\sqrt{63+1}$.
\end{remark}

\begin{conjecture}\label{Conj2.20} {\rm  (2013-10-15) (i)} For any integer $n>5$, there
is a prime $p<n$ with $p+6$ and $n+(n-p)^2$ both prime.

{\rm (ii)} For any integer $n>3$, there is a prime $p<n$ with $3p-4$ and
$n^2+(n-p)^2$ both prime.
\end{conjecture}
\begin{remark}\label{Rem2.20} See \cite[A227898 and A227899]{S} for related data.
\end{remark}

\begin{conjecture}\label{Conj2.21} {\rm (2013-11-20) (i)} Any integer $n>1$ can be written
as $x+y$ with $x,y\in\Z^+$ such that $x+ny$ and $x^2+ny^2$ are both
prime.

{\rm (ii)} Any integer $n>2$ can be written
as $x+y$ with $x,y\in\Z^+$ such that $nx+y$ and $nx-y$ are both prime.
Also, any integer $n>2$ can be written
as $x+y$ with $x,y\in\Z^+$ such that $x^2+(n-2)y^2$ is prime.

{\rm (iii)} Any integer $n>2$ can be written
as $p+q$ with $q\in\Z^+$ such that $p$ and $p^3+nq^2\ (or\ p+nq)$ are both prime.
\end{conjecture}
\begin{remark}\label{Rem2.21} See \cite[A232174, A231883 and A232186]{S} for related data. For example, $20 = 11+9$ with
$11+20\cdot9=191$ and $11^2+20\cdot9^2=121+20\times81=1741$ both prime.
The author would like to offer 200 US dollars as the prize for the first solution to part (i) of Conjecture \ref{Conj2.21}.
We also conjecture that there are infinitely many $n\in\Z^+$
such that $p_n=x^2+ny^2$ for some $x,y\in\Z^+$ (where $p_n$ is the $n$-th prime).
\end{remark}

\begin{conjecture}\label{Conj2.22} {\rm  (2013-10-14)} Any integer $n>1$ can be written as $x+y$ with $x,y\in\Z^+$
such that $2x+1$, $x^2+x+1$ and $y^2+y+1$ are all prime. Also, each integer $n>1$ can be written as $x+y$ with $x,y\in\Z^+$
such that $x^2+1\ (or\ 4x^2+1)$ and $4y^2+1$ are both prime.
\end{conjecture}
\begin{remark}\label{Rem2.22} See \cite[A230252]{S} for related data. For example, $31=14+17$ with $2\cdot14+1=29$, $14^2+14+1=211$
and $17^2+17+1=307$ all prime.
\end{remark}

In 2001 Heath-Brown \cite{HB} proved that there are infinitely many primes
of the form $x^3+2y^3$ where $x$ and $y$ are positive integers.

\begin{conjecture}\label{Conj2.23} {\rm  (2012-12-14)} Any positive integer $n$ can be
written as $x+y\ (x,y\in\N)$ with $x^3+2y^3$ prime.
In general, for each positive {\it odd} integer $m$, any
sufficiently large integer can be written as $x+y\ (x,y\in\N)$ with
$x^m+2y^m$ prime.
\end{conjecture}
\begin{remark}\label{Rem2.23} See \cite[A220413]{S} for related data. For any integer $d>2$, not every sufficiently large integer $n$
can be written as $x+y\ (x,y\in\N)$ with $x^3+dy^3$ prime. For, if $n$
is a multiple of a prime divisor $p$ of $d-1$, then
$x^3+d(n-x)^3\eq (1-d)x^3\eq0\pmod{p}$ for any integer $x$.
\end{remark}

\begin{conjecture}\label{Conj2.24} {\rm  (2013-04-15)} For any integer $n>4$, there is a
positive integer $k<n$ such that $p=2n+k$ and
$2n^3+k^3=2n^3+(p-2n)^3$ are both prime.
\end{conjecture}
\begin{remark}\label{Rem2.24} See \cite[A224030]{S} for related data.
\end{remark}

\begin{conjecture}\label{Conj2.25} {\rm (2012-12-16)} Let $m$ be a positive integer.
Then, any sufficiently large odd integer $n$ can be written as $x+y\
(x,y\in\Z^+)$ with $x^m+3y^m$ prime (and any sufficiently large even
integer $n$ can be written as $x+y\ (x,y\in\Z^+)$ with $x^m+3y^m+1$
prime). In particular, if $m\ls 6$ or $m=18$, then  each
positive odd integer can be written as $x+y\ (x,y\in\N)$ with
$x^m+3y^m$ prime.
\end{conjecture}
\begin{remark}\label{Rem2.25} See \cite[A220572]{S} for related data and comments. For example, $5$ can be written as $1+4$ with
$$1^{18}+3\cdot 4^{18}=206158430209$$ prime.
\end{remark}

\section{Other representation problems for positive integers}
\setcounter{lemma}{0}
\setcounter{theorem}{0}
\setcounter{corollary}{0}
\setcounter{remark}{0}
\setcounter{equation}{0}
\setcounter{conjecture}{0}

\begin{conjecture}\label{Conj3.1}  {\rm (2013-11-12)} Any integer $n>1$ can be written
as $x+y\ (x,y\in\Z^+)$ with $[x,y]+1$ prime, where $[x,y]$ is the least common multiple of $x$ and $y$.
Also, each integer $n>3$ can be written as $x+y\ (x,y\in\Z^+)$ with $[x,y]-1$ prime.
\end{conjecture}
\begin{remark}\label{Rem3.1} See \cite[A231635]{S} for related data. For example, $10=4+6$ with $[4,6]+1=13$
and $[4,6]-1=11$ both prime.
\end{remark}

As usual, for $x\in\Z$ we let $T_x$ denote the triangular number $x(x+1)/2$.

\begin{conjecture}\label{Conj3.2}  {\rm (i) (2013-11-10)} Any integer $n>1$ can be written
as $x+y\ (x,y\in\Z^+)$ with $T_x+y^2$ prime. Also, any integer $n>6$ can be written
as $x+y\ (x,y\in\Z^+)$ with $T_x+y^4$ prime.

{\rm (ii) (2013-11-18)} Any integer $n>1$ can be written
as $x+y\ (x,y\in\Z^+)$ with $p=2x+1$ and $T_x+y=n+(p-1)(p-3)/8$ both prime.
\end{conjecture}
\begin{remark}\label{Rem3.2} See \cite[A228425 and A232109]{S} for related data.
For example, $18=7+11$ with $T_7+11^2=149$ prime, $27=5+22$ with $T_5+22^4=234271$ prime,
and $18=11+7$ with $2\cdot11+1=23$
and $T_{11}+7=73$ both prime.
\end{remark}

\begin{conjecture}\label{Conj3.3}  {\rm (2012-10-15)} Each $n=1,2,3.\ldots$ can be written
as $T_x+y$ with $x,y\in\N$ such that $T_y+1$ is prime.
\end{conjecture}
\begin{remark}\label{Rem3.3} See \cite[A229166]{S} for related data.
For example, $34$ has a unique required representation: $34=T_5+19$ with $T_{19}+1=191$ prime.
\end{remark}

\begin{conjecture}\label{Conj3.4}  {\rm (2012-12-09)} Any integer $n>2$ can be written
as $x^2+y\ (x,y\in\Z^+)$ with $2xy-1$ prime. In other words, for
each $n=3,4,\ldots$ there is a prime of the form $2k(n-k^2)-1$ with
$k\in\Z^+$.
\end{conjecture}
\begin{remark}\label{Rem3.4} We have verified Conjecture \ref{Conj3.4} for all $n=3,4,\ldots,3\cdot10^9$. See \cite[A220272]{S} for related data.
For example, $18=3^2+9$ with $2\times3\times9-1=53$ prime.
\end{remark}

\begin{conjecture}\label{Conj3.5}  {\rm (2013-10-21)} Any integer $n>1$ can be written
as $x^2+y$ with $2y^2-1$ prime, where $x,y\in\N$. In other words, for each $n=2,3,4,\ldots$ there is an integer $0\ls k\ls\sqrt n$
such that $2(n-k^2)^2-1$ is prime.
\end{conjecture}
\begin{remark}\label{Rem3.5} We have verified this conjecture for all $n=2,3,\ldots,10^8$.
See \cite[A230494]{S} for related data. For example, $9=1^2+8$ with $2\cdot8^2-1=127$ prime.
\end{remark}

\begin{conjecture}\label{Conj3.6}  {\rm (2013-11-11) (i)} Any integer $n>1$ can be written
as $k+m\ (k,m\in\Z^+)$ with $2^k+m$ prime. In other words, for each
$n=2,3,\ldots$ there is a positive integer $k<n$ with $n+2^k-k$
prime.

{\rm (ii)} For any integer $n>3$, there is a positive integer $k<n$ such that $n+2^k+k$ is prime.
\end{conjecture}
\begin{remark}\label{Rem3.6} We have verified parts (i) and (ii) of this conjecture for $n$ up to $10^7$ and $3.8\times10^6$ respectively,
see \cite[A231201, A231557 and A231725]{S} for related data and other similar conjectures. For example, $9302003=311468+8990535$
with $2^{311468}+8990535$ a prime of 93762 decimal digits.
In \cite{S13c} the author proved that the set $\{2^k-k: k=1,2,3,\ldots\}$ contains a complete system of residues modulo any positive integer.
The author would like to offer 1000 US dollars as the prize for the first solution to part (i) of Conjecture \ref{Conj3.6}.
\end{remark}

\begin{conjecture}\label{Conj3.7}  {\rm (2013-11-23)} Any integer $n>3$ can be written as $p+(2^k-k)+(2^m-m)$ with $p$ prime and $k,m\in\Z^+$.
\end{conjecture}
\begin{remark}\label{Rem3.7}  For example, $94$ has a unique required representation $31+(2^3-3)+(2^6-6)$. See \cite[A232398]{S} for related data.
After the author verified this conjecture for $n$ up to $2\times10^8$, Qing-Hu Hou extended the verification to $10^{10}$ in Dec. 2013.
In contrast with Conjecture \ref{Conj3.7}, R. Crocker \cite{Cr} proved in 1971 that there are infinitely many positive odd numbers not of the form $p + 2^k + 2^m$
with $p$ prime and $k,m\in\Z^+$.
\end{remark}

\begin{conjecture}\label{Conj3.8}  {\rm (2013-11-11)} Let $r\in\{1,2\}$. Then any integer $n>1$ can be written
as $k+m\ (k,m\in\Z^+)$ with $2^km^r+1$ prime. Also, any integer $n>2$ can be written
as $k+m\ (k,m\in\Z^+)$ with $2^km^r-1$ prime.
\end{conjecture}
\begin{remark}\label{Rem3.8} See \cite[A231561]{S} for related data.
\end{remark}

\begin{conjecture}\label{Conj3.9} {\rm (i) (2013-11-10)} Any integer $n>1$ can be written
as $k+m\ (k,m\in\Z^+)$ with $F_k+m$ $($or $F_k+2m$, or $F_k+m(m+1))$ prime, where the Fibonacci sequence $(F_j)_{j\gs0}$ is given by $F_0=0$, $F_1=1$, and $F_{j+1}=F_j+F_{j-1}$
for $j\in\Z^+$.

 {\rm (ii) (2014-04-27)} Any integer $n>1$ can be written
as $k+m\ (k,m\in\Z^+)$ with $L_k+m$ prime, where the Lucas sequence $(L_j)_{j\gs0}$ is given by $L_0=2$, $L_1=1$, and $L_{j+1}=L_j+L_{j-1}$
for $j\in\Z^+$.
\end{conjecture}
\begin{remark}\label{Rem3.9} See \cite[A231555 and A241844]{S} for related data.
We have verified parts (i) and (ii) of Conjecture \ref{Conj3.9} for $n$ up to $3.7\times10^6$ and $7\times10^6$ respectively.
\end{remark}

\begin{conjecture}\label{Conj3.10}  {\rm (i) (2013-11-11)} Any integer $n>1$ can be written as $k+m\ (k,m\in\Z^+)$ with $k!m+1$ prime. Also,
Any integer $n>3$ can be written as $k+m\ (k,m\in\Z^+)$ with $k!m-1$ prime.

{\rm (ii) (2014-03-19)} Let $r\in\{1,-1\}$. For each integer $n>1$, there is a number $k\in\{1,\ldots,n\}$ with $k!n+r$ prime.
\end{conjecture}
\begin{remark}\label{Rem3.10} See \cite[A231516 and A231631]{S} for related data. We have verified part (i) of Conjecture \ref{Conj3.10} for $n$ up to $10^6$.
We also conjecture that for any integer $n>2$ there is a positive integer $k<\sqrt{n}\log n$ with $k!(n-k)+1$ prime.
\end{remark}

\begin{conjecture}\label{Conj3.11} {\rm (i) (2015-04-01)} Let $k,m\in\Z^+$ with $k+m>2$. Then any integer $n>2$ can be written as $\lfloor p/k\rfloor+\lfloor q/m\rfloor$
with $p$ and $q$ both prime.

{\rm (ii) (2015-04-24)} Let
\begin{align*} T:=&\l\{\l\lfloor\f{x}9\r\rfloor:\ x-1\ \t{and}\ x+1\ \t{are twin prime}\r\}
\\=&\l\{\l\lfloor\f{x}3\r\rfloor:\ 3x-1\ \t{and}\ 3x+1\ \t{are twin prime}\r\}.
\end{align*}
Then, any positive integer can be written as the sum of two distinct elements of $T$ one of which is even.
\end{conjecture}
\begin{remark}\label{Rem3.11} See \cite[A256555 and A256707]{S} for related data. Part (i) of Conjecture \ref{Conj3.11} in the case $k=m=2$ reduces to Goldbach's conjecture,
and it reduces to Lemoine's conjecture when $\{k,m\}=\{1,2\}$. Part (ii) of Conjecture \ref{Conj3.11} implies the twin prime conjecture.
\end{remark}

\begin{conjecture}\label{Conj3.12} {\rm (2014-03-03)} {\rm (i)} Let $1<m<n$ be integers with $m\nmid n$. Then $\lfloor kn/m\rfloor$ is prime for some $k=1,\ldots,n-1$.

{\rm (ii)} Let $m>2$ and $n>2$ be integers. Then there is a prime $p<n$ with $\lfloor(n-p)/m\rfloor$ a square.
Also, there is a prime $p<n$ such that $\lfloor(n-p)/m\rfloor$ is a triangular number of the form $T_{(q-3)/2}=(q-1)(q-3)/8$ with $q$ an odd prime.

{\rm (iii)} For each $n=3,4,\ldots$, there is a prime $p<n$ with $\lfloor(n-p)/5\rfloor$ a cube.
\end{conjecture}
\begin{remark}\label{Rem3.12} See \cite[A238703, A238732 and A238733]{S} for related data.
\end{remark}

\begin{conjecture}\label{Conj3.13}  {\rm (i) (2013-10-21)} Let
$$S=\{n\in\Z^+:\ 2n+1\ \mbox{and}\ 2n^3+1\ \mbox{are both prime}\}.$$
Then any integer $n>2$ is a sum of three elements of $S$.

{\rm (ii) (2013-10-22)} Any integer $n>5$ can be written as $a+b+c$ with $a,b,c\in\Z^+$ such that
$$\{a^2+a\pm1\},\ \{b^2+b\pm1\},\ \{c^2+c\pm1\}$$
are all twin prime pairs!
\end{conjecture}
\begin{remark}\label{Rem3.13} See \cite[A230507 and A230516]{S} for related data and comments.
\end{remark}

\begin{conjecture}\label{Conj3.14} {\rm (2013-10-11)} Let
$$P=\{p:\ p,\ p+6\ \mbox{and}\ 3p+8\ \mbox{are all prime}\}.$$
Then, for any integer $n>6$, we can write $2n+1=p+q+r$ with
$p,q,r\in P$ such that $p+q+9$ is also prime.
\end{conjecture}
\begin{remark}\label{Rem3.14} This implies not only Goldbach's weak conjecture but
also Goldbach's conjecture for even numbers. See \cite[A230217 and A230219]{S} for related data.
Note that $37$ has a unique required representation $7+13+17$; in fact,
\begin{align*}&7,\ 7+6=13,\ 3\times7+8=29,
\\&13,\ 13+6=19,\ 3\cdot 13+8=47,
\\&17,\ 17+6=23,\ 3\cdot 17+8=59,
\end{align*}
and $7+13+9 = 29$ are all prime.
\end{remark}

\begin{conjecture}\label{Conj3.15} {\rm (2013-10-12)} Let
$$P'=\{p:\ p,\ 3p-4,\ 3p-10\ \mbox{and}\ 3p-14\ \mbox{are all prime}\}.$$
Then, for any integer $n>17$, we can write $2n=p+q+r+s$ with
$p,q,r,s\in P'$.
\end{conjecture}
\begin{remark}\label{Rem3.15} See \cite[A230223 and A230224]{S} for related data.
Note that such a representation involves 16 primes!
For example, $54$ has a unique required representation $7+11+17+19$; in fact,
\begin{align*}&7,\ 3\cdot7-4=17,\ 3\cdot7-10=11,\ 3\cdot7-14=7,
\\&11,\ 3\cdot11-4=29,\ 3\cdot11-10=23,\ 3\cdot11-14=19,
\\&17,\ 3\cdot17-4=47,\ 3\cdot17-10=41,\ 3\cdot17-14=37,
\\&19,\ 3\cdot19-4=53,\ 3\cdot19-10=47,\ 3\cdot19-14=43
\end{align*}
are all prime.
\end{remark}

\begin{conjecture}\label{Conj3.16} {\rm (i) (2015-10-01)} Any integer $n>1$ can be written as $x^2+y^2+\varphi(z^2)$ with $x,y\in\N$, $x\ls y$ and $z\in\Z^+$
such that $y$ or $z$ is prime.

{\rm (ii) (2015-10-02)} Each positive integer can be written as $x^2+y^2+p(p+\ve)/2$, where $x,y\in\Z$, $\ve\in\{\pm1\}$, and $p$ is a prime.
\end{conjecture}
\begin{remark}\label{Rem3.16} See \cite[A262311, A262785, A262982, A262985, A263992, A263998, A264010 and A264025]{S} for related data and similar conjectures.
For example, $13=1^2+2^2+\varphi(4^2)$ with $2$ prime, $94415=115^2+178^2+\varphi(223^2)$ with $223$ prime, $97=1^2+9^2+5(5+1)/2$ with $5$ prime, and $538=3^2+8^2+31(31-1)/2$
with $31$ prime. It is known that each $n\in\N$ can be expressed as the sum of two squares and a triangular number (cf. \cite{S07}).
\end{remark}

\begin{conjecture}\label{Conj3.17} {\rm (2014-02-26) (i)} Any integer $n>6$ can be written as $k+m\ (k,m\in\Z^+)$ such that $p_k+\pi(m)$ is a triangular number.

{\rm (ii)} Any integer $n>10$ can be written as $k+m\ (k,m\in\Z^+)$ such that $p=p_k+\pi(m)$ and $p+2$ are both prime.
\end{conjecture}
\begin{remark}\label{Rem3.17} See \cite[A238405 and A238386]{S} for related data. For example, $72=41+31$ with $p_{72}+\pi(31)=179+11=19\cdot{20}/2$ a triangular number, and
$108=15+93$ with $p_{15}+\pi(93)=47+24=71$ and $71+2=73$ twin prime.
\end{remark}

\begin{conjecture}\label{Conj3.18} {\rm (2014-03-05)} Any integer $n>2$ can be written as $q+m$ with $m\in\Z^+$ such that
$q$, $p_q-q+1$ and $p_{p_m}-p_m+1$ are all prime.
\end{conjecture}
\begin{remark}\label{Rem3.18} See \cite[A237715]{S} for related data.
\end{remark}

\begin{conjecture}\label{Conj3.19} {\rm (2014-01-04)} For any integer $n>6$, there is a prime $q<n/2$ with $p_q-q+1$ prime
 such that $n-(1+\{n\}_2)q$ is prime.
\end{conjecture}
\begin{remark}\label{Rem3.19} This conjecture is stronger than Goldbach's conjecture and Lemoine's conjecture, and it also implies that
there are infinitely many primes $q$ with $p_q-q+1$ prime.
See \cite[A235189]{S} for related data. For example, $7$, $p_7-7+1=17-6=11$ and $61-2\cdot7=47$ are all prime,
and $31$, $p_{31}-31+1=97$ and $98-31=67$ are all prime.
\end{remark}

\begin{conjecture}\label{Conj3.20} {\rm (2014-02-04) (i)} For any integer $n>2$, we can write $2n=p+q$ with $p$, $q$ and $\varphi(p+2)\pm1$ all prime.
Also, for any integer $n>12$ we can write $2n-1=2p+q$
with $p$, $q$ and $\varphi(p+1)\pm1$ all prime.

{\rm (ii)} Any integer $n\gs24$ can be written as $(1+\{n\}_2)p+q$ with $p,q,\varphi(p+1)-1$ and $\varphi(q-1)+1$ all prime.
\end{conjecture}
\begin{remark}\label{Rem3.20}
See \cite[A237168, A237183 and A237184]{S} for related data. Note that either of the two parts is stronger than Goldbach's conjecture and Lemoine's conjecture.
Also, part (i) of Conjecture \ref{Conj3.20} implies the twin prime conjecture.
\end{remark}

\begin{conjecture}\label{Conj3.21} {\rm (2014-02-04) (i)} Any integer $n\gs12$ can be written as $k+m$ with
$k,m\in\Z^+$ and $k\not=m$ such that $\varphi(k)\pm1$ and $\varphi(m)\pm1$ are all prime.

{\rm (ii)} Any integer $n>6$ can be written as $k+m\ (k,m\in\Z^+)$
such that both $\{p_k,p_k+2\}$ and $\{\varphi(m)-1,\varphi(m)+1\}$ are twin prime pairs.

{\rm (iii)} Any integer $n\gs6$ can be written as $k+m\ (k,m\in\Z^+)$ with $p_{p_k}+2$ and $\varphi(m)\pm1$ all prime. Also,
each $n=2,3,4,\ldots$ can be written as $k+m\ (k,m\in\Z^+)$ with $p_{p_k}+2$ and $6m\pm1$ all prime.

{\rm (iv)} Any integer $n\gs8$ can be written as $k+m\ (k,m\in\Z^+)$ with $p_{p_{p_k}}-2$ and $\varphi(m)\pm1$ all prime.

{\rm (v)} Any integer $n>8$ can be written as $k+m\ (k,m\in\Z^+)$ with $3k\pm1$ and $\varphi(m)\pm1$ all prime.

{\rm (vi)} Any integer $n\gs12$ can be written as $p+q\ (q\in\Z^+)$ with $p$, $p+6$ and $\varphi(q)\pm1$ all prime.
\end{conjecture}
\begin{remark}\label{Rem3.21}
See \cite[A237127, A237130, A218829 and A237253]{S} for related data and comments. Clearly each part of Conjecture \ref{Conj3.21} implies the twin prime conjecture.
\end{remark}

\begin{conjecture}\label{Conj3.22} {\rm (2013-12-31) (i)} Any integer $n>1$ with $n\not=8$ can be written as $k+m\ (k,m\in\Z^+)$ such that $p=k+\varphi(m)$ and $2n-p$
are both prime.

{\rm (ii)} Each integer $n>2$ can be written as $k+m\ (k,m\in\Z^+)$ such that $p=k+\varphi(m)$ and $2n+1-2p$
are both prime.
\end{conjecture}
\begin{remark}\label{Rem3.22} Clearly parts (i) and (ii) are stronger than Goldbach's conjecture and Lemoine's conjecture respectively.
See \cite[A234808 and A234809]{S} for related data. For example, $24=9+15$ with $9+\varphi(15)=17$ and $2\cdot24-17=31$ both prime,
and $41=7+34$ with $7+\varphi(34)=23$ and $2\cdot41+1-2\cdot23=37$ both prime.
\end{remark}

\begin{conjecture}\label{Conj3.23} {\rm (2014-02-02) (i)} Any integer $n>1$ can be written as $k+m\ (k,m\in\Z^+)$ such that $6k\pm1$ and $k+\varphi(m)$ are all prime.

{\rm (ii)} Any integer $n>3$ with $n\not=12$ can be written as $k+m\ (k,m\in\Z^+)$ such that $6k\pm1$ and $k+\varphi(m)/2$ are all prime.

{\rm (iii)} Each integer $n>5$ can be written as $k+m\ (k,m\in\Z^+)$ with $k+\varphi(m)/2$ a square.
\end{conjecture}
\begin{remark}\label{Rem3.23} See \cite[A236968 and A236567]{S} for related data.
\end{remark}

\begin{conjecture}\label{Conj3.24} {\rm (2014-01-13)} Define
$$K:=\{k\in\Z^+:\ k(k+1)-p_k\ \mbox{is prime}\}.$$

{\rm (i)} Any integer $n>3$ can be written as $a+b$ with $a,b\in K$.

{\rm (ii)} Any integer $n>2$ can be expressed as the sum of an element of $K$ and a positive triangular number.

{\rm (iii)} Any integer $n>3$ can be written as the sum of an element of $K$ and a prime $q$ with $p_q-q+1$ also prime.

{\rm (iv)} Any integer $n>7$ can be written as $k+m\ (k,m\in\Z^+)$ such that $q=p_k+\varphi(m)$ and $q(q+1)-p_q$ are both prime.
\end{conjecture}
\begin{remark}\label{Rem3.24} See \cite[A235592, A235613, A235614, A235661, A235703, A232353]{S} for related data.
\end{remark}

\begin{conjecture}\label{Conj3.25} {\rm (2014-01-18) (i)} Any integer $n>7$ can be written as $k+m\ (k,m\in\Z^+)$ such that $p=\varphi(k)+\varphi(m)/2-1$
is a prime and also $2$ is a primitive root modulo $p$.

{\rm (ii)} Any integer $n\gs38$ can be written as $k+m\ (k,m\in\Z^+)$ such that $p=p_k+\varphi(m)$
is a Sophie Germain prime and also $2$ is a primitive root modulo $p$.
\end{conjecture}
\begin{remark}\label{Rem3.25} See \cite[A235987]{S} for related data and comments. For example,
$79=19+60$, $p_{19}+\varphi(60)=67+16=83$ is a Sophie Germain prime and $2$ is a primitive root modulo $83$.
\end{remark}

\begin{conjecture}\label{Conj3.26} {\rm (i) (2012-12-23)} Any integer $n>5$ can be written as $k+m\ (k,m\in\{3,4,\ldots\})$ with $2^{\varphi(k)}+2^{\varphi(m)}-1$ prime.

{\rm (ii) (2012-12-24)} For any integer $a>1$, there is a positive integer $N(a)$ such that
any integer $n > N(a)$ can be written as $k+m$ with $k,m\in\{3,4,\ldots\}$
 such that $a^{\varphi(k)} + a^{\varphi(m)/2} - 1$ is prime. Moreover, we may take $N(2) = N(3) = \ldots = N(6) = N(8) = 5$ and $N(7) = 17$.
\end{conjecture}
\begin{remark}\label{Rem3.26} See \cite[A234309, A234347 and A234359]{S} for related data and comments.
Clearly, part (ii) of Conjecture \ref{Conj3.26} implies that for each $a =2,3,\ldots$ there are infinitely many primes of the form $a^{2k}+a^m-1$ with $k,m\in\Z^+$.
\end{remark}

\begin{conjecture}\label{Conj3.27} {\rm (2013-12-26) (i)} Any integer $n\gs10$ can be written as $k+m\ (k,m\in\Z^+)$ with $2^{\varphi(k)/2+\varphi(m)/6}+3$
prime. Also, any integer $n>13$ can be written as $k+m\ (k,m\in\Z^+)$ with $2^{\varphi(k)/2+\varphi(m)/6}-3$ prime.

{\rm (ii)} Any integer $n>25$ can be written as $k+m\ (k,m\in\Z^+)$ with $3\cdot2^{\varphi(k)/2+\varphi(m)/8}+1$
prime. Also, any integer $n\gs15$ can be written as $k+m\ (k,m\in\Z^+)$ with $3\cdot2^{\varphi(k)/2+\varphi(m)/12}-1$
prime.

{\rm (iii)} Any integer $n\gs27$ can be written as $k+m\ (k,m\in\Z^+)$ with $2\cdot3^{\varphi(k)/2+\varphi(m)/12}$
$+1$ prime. Also, any integer $n>37$ can be written as $k+m\ (k,m\in\Z^+)$ with $2\cdot3^{\varphi(k)/2+\varphi(m)/12}-1$
prime.

{\rm (iv)} Any integer $n>10$ can be written as $k+m\ (k,m\in\Z^+)$ with $2^{\varphi(k)+\varphi(m)/4}-5$ prime.
\end{conjecture}
\begin{remark}\label{Rem3.27} This implies that there are infinitely many primes in any of the following seven forms:
$$2^n+3,\ 2^n-3,\ 3\cdot2^n+1,\ 3\cdot2^n-1,\ 2\cdot3^n+1,\ 2\cdot3^n-1,\ 2^n-5.$$
We have verified Conjecture \ref{Conj3.27} for $n$ up to 50,000. See \cite[A234451, A236358 and A234504]{S} for related data.
\end{remark}

\begin{conjecture}\label{Conj3.28} {\rm (2012-12-24) (i)} Any integer $n>1$ can be written as $k+m\ (k,m\in\Z^+)$ with $(k+1)^{\varphi(m)}+k$ prime.
Also, each integer $n>1$ can be written as $k+m\ (k,m\in\Z^+)$ with $k(k+1)^{\varphi(m)}+1$ prime.

{\rm (ii)} Any integer $n>5$ can be written as $k+m\ (k,m\in\Z^+)$ with $(k+1)^{\varphi(m)/2}-k$ prime.
Also, each integer $n>3$ can be written as $k+m\ (k,m\in\Z^+)$ with $k(k+1)^{\varphi(m)/2}-1$ prime.
\end{conjecture}
\begin{remark}\label{Rem3.28} This conjecture is somewhat curious. See \cite[A234360]{S} for related data.
\end{remark}

\begin{conjecture}\label{Conj3.29} {\rm (i) (2014-02-02)} Any integer $n>8$ can be written as $i+j$ with $i,j\in\Z^+$ and $i<j$ such that $\varphi(i)\varphi(j)$ is a square.
Also, for each $k=3,4,\ldots$, any integer $n\gs 3k$ can be written as $i_1+i_2+\ldots+i_k$ with $i_1,i_2,\ldots,i_k\in\Z^+$ not all equal such that $\varphi(i_1)\varphi(i_2)\cdots\varphi(i_k)$ is a $k$-th power.

{\rm (ii) (2014-02-09)} Any integer $n\gs8$ can be written as $i+j$ with $i,j\in\Z^+$ and $i<j$ such that $\varphi(ij)+1$ is a square.
Also, for each $k=3,4,\ldots$, any integer $n>2k+1$ can be written as $\sum_{j=1}^ki_j$ with $i_1,i_2,\ldots,i_k\in\Z^+$ such that $\varphi(i_1i_2\ldots i_k)$
is a $k$-th power.

{\rm (iii) (2014-02-04)} Let $k>1$ be an integer. Any sufficiently large integer $n$ can be written as $\sum_{j=1}^ki_j$ with $i_1,\ldots,i_k\in\Z^+$ and $i_1<\ldots<i_k$
such that all those $\varphi(i_j)\ (j=1,\ldots,k)$ are $k$-th powers.

{\rm (iv) (2014-02-02)} For each $k=3,4,\ldots$ any sufficiently large integer $n$ can be written as $i_1+i_2+\ldots+i_k$ with $i_1,i_2,\ldots,i_k$ not all equal such that
the product $i_1i_2\ldots i_k$ is a $k$-th power.
\end{conjecture}
\begin{remark}\label{Rem3.29} See \cite[A236998, A233386, A237523, A237524, A237123, A237050]{S} for related data.
For any integer $k>1$, we clearly have $2k+2=4+(k-1)2$ with $\varphi(4\cdot2^{k-1})=2^k$ a $k$-th power.
In contrast with part (i) of Conjecture \ref{Conj3.29}, we also conjecture that (cf. \cite[A237049]{S}) for each $k=2,3,4,\ldots$ any sufficiently large integer $n$ can be written as
$\sum_{j=1}^ki_j$ with $i_1,i_2,\ldots,i_k\in\Z^+$ not all equal such that $\prod_{j=1}^k\sigma(i_j)$ is a $k$-th power, where $\sigma(m)$ denotes the sum of all positive divisors of $m\in\Z^+$.
\end{remark}

\begin{conjecture}\label{Conj3.30} {\rm (i) (2013-12-21)} Any integer $n>5$ can be written as $k+m$ with $k,m\in\Z^+$ such that $(\varphi(k)+\varphi(m))/2$ is prime.

{\rm (ii) (2013-12-22)} Any positive integer $n$ not dividing $6$ can be written as $k+m$ with $k,m\in\Z^+$ such that $k\varphi(m)+1$ is a square.
Also, any integer $n>4$ can be written as $k+m$ with $k,m\in\Z^+$ and $k<m$ such that $k\varphi(m)-1$ and $k\varphi(m)+1$ are both prime.

{\rm (iii) (2013-12-12)} Any integer $n>5$ can be written as $k+m$ with $k,m\in\Z^+$ such that $\varphi(k)\varphi(m)-1$ and $\varphi(k)\varphi(m)+1$ are both prime.

{\rm (iv) (2013-12-23)} Any integer $n>4$ can be written as $k+m\ (k,m\in\Z^+)$ with $\varphi(k^2)\varphi(m)-1$ a Sophie Germain prime.
\end{conjecture}
\begin{remark}\label{Rem3.30} See \cite[A233918, A234200, A234246, A233547, A234308]{S} for related data. For example, $13=3+10$ with $(\varphi(3)+\varphi(10))/2 = 3$ prime,
$13=4+9$ with $4\varphi(9)+1=25$ a square, $18=5+13$ with $\{5\varphi(13)\pm1\}=\{59,61\}$ a twin prime pair,
$26=7+19$ with $\{\varphi(7)\varphi(19)\pm1\}=\{107,109\}$ a twin prime pair, and $30=2+28$ with $\varphi(2^2)\varphi(28)-1=23$ a Sophie Germain prime.
\end{remark}

\begin{conjecture}\label{Conj3.31} {\rm (2013-12-12) (i)} Any integer $n>1$ can be written as $k^2+m$ with $\sigma(k^2)+\varphi(m)$ prime, where
$k$ and $m$ are positive integers with $k^2\ls m$.

{\rm (ii)} Any integer $n>1$ can be written as $k+m$ with $k,m\in\Z^+$ such that $\sigma(k)^2+\varphi(m)\ ($or $\sigma(k)+\varphi(m)^2)$ is prime.
\end{conjecture}
\begin{remark}\label{Rem3.31} See \cite[A233544]{S} for related data and comments.  We have verified part (i) of Conjecture \ref{Conj3.31} for all $n=2,\ldots,10^8$;
for example, $25=2^2+21$ with $\sigma(2^2)+\varphi(21)=7+12=19$ prime.
\end{remark}

\begin{conjecture}\label{Conj3.32} Let $n>2$ be an integer.

{\rm (i) (2013-12-14)} If $n$ is even, then $n$ can be written as $p+\sigma(k)$, where $p$ is an odd prime and $k\in\{1,\ldots,n-1\}$.

{\rm (ii) (2013-12-17)} If $n$ is odd, then $n$ can be written as $p+\varphi(k^2)$, where $p$ is a prime and $k$ is a positive integer with $k^2<n$.
\end{conjecture}
\begin{remark}\label{Rem3.32} See \cite[A233654, A233793 and A233867]{S} for related data. For example, $28=13+\sigma(8)$ with $13$ prime, and $29=23+\varphi(3^2)$ with $23$ prime.
Note that if $n=p+q$ with $p$ and $q$ both prime then $n+1=p+(q+1)=p+\sigma(q)$ and $n-1=p+(q-1)=p+\varphi(q)$.
\end{remark}

\begin{conjecture}\label{Conj3.33} {\rm (2012-12-29) (i)} For each integer $n>8$ with $n\not=14$, there is a prime $p$ between $n$ and $2n$
with $(\f np)=1$. If $n\in\Z^+$ is not a square, then there is a prime $p$ between $n$ and $2n$ with $(\f np)=-1$.

 {\rm (ii)} For any integer $n>5$ there is a prime $p\in(n,2n)$ with $(\f{2n}p)=1$.
 For any integer $n>6$ there is a prime $p\in(n,2n)$ with $(\f{-n}p)=-1$.
\end{conjecture}
\begin{remark}\label{Rem3.33} We have verified this refinement of Bertrand's postulate for $n$ up to $5\times 10^8$.
\end{remark}

\begin{conjecture}\label{Conj3.34} {\rm (2012-12-29)} For any positive integer $n$ there is a prime $p$ between $n^2$ and $(n+1)^2$
with $(\f np)=1$.
Also, for any integer $n>1$ we have $(\f {n(n+1)}p)=1$ for some prime $p\in(n^2,(n+1)^2)$.
\end{conjecture}
\begin{remark}\label{Rem3.34} We have verified this refinement of Legendre's conjecture for $n$ up to $10^9$.
\end{remark}

\begin{conjecture}\label{Conj3.35} {\rm (Olivier Gerard and Zhi-Wei Sun, 2012-11-19)}  For any integer $n\gs400$ with $n\not=757,\, 1069,\, 1238$,
there are odd primes $p$ and $q$ with $(\f pq)=(\f qp)=1$ such that $p+(1+\{n\}_2)q=n$.
\end{conjecture}
\begin{remark}\label{Rem3.35}  We have verified Conjecture \ref{Conj3.35} for $n$ up to $10^8$. See \cite{GS} for the announcement of this conjecture.
\end{remark}

\begin{conjecture}\label{Conj3.36} {\rm (2012-11-22)} Let $m$ be any integer. Then, for every sufficiently large integer $n$
there are primes $p>q>2$ with $(\f{p-(1+\{n\}_2)m}q)=(\f{q+m}p)=1$ and $p+(1+\{n\}_2)q=n$.
\end{conjecture}
\begin{remark}\label{Rem3.36}  Conjecture \ref{Conj3.36} in the case $m=0$ corresponds to Conjecture \ref{Conj3.35}.
\end{remark}

\begin{conjecture}\label{Conj3.37} {\rm (2012-12-30)}  Any integer $n>5$ can be written as $p+(1+\{n\}_2)q$, where $p$ is an odd prime
and $q$ is a prime not exceeding $n/2$ such that $(\f qn)=1$ if $2\nmid n$,
and $(\f{(q+1)/2}{n+1})=1$ if $2\mid n$.
\end{conjecture}
\begin{remark}\label{Rem3.37} We have verified this refinement of Goldbach's and Lemoine's conjectures for $n$ up to $10^9$.
\end{remark}

\begin{conjecture}\label{Conj3.38} {\rm (2013-01-19) (i)} Any even integer $2n>4$ can be written as $p+q=(p+1)+(q-1)$, where $p$ and $q$ are primes with $p+1$ and $q-1$
both practical.

{\rm (ii)} For each integer $n>8$, we can write $2n-1=p+q=2p+(q-p)$, where $p$ and $q-p$ are both prime, and $q$ is practical.
\end{conjecture}
\begin{remark}\label{Rem3.38} We have verified both parts of Conjecture \ref{Conj3.38} for $n$ up to $10^8$. See \cite[A209320 and A209315]{S}
for related data.
\end{remark}

If one of $n$ and $n+1$ is prime and the other is practical, then we call $\{n,n+1\}$ a {\it couple}.
As powers of two are practical numbers, $\{2^p-1,2^p\}$ is a couple if $2^p-1$ is a Mersenne prime, and $\{2^{2^n},2^{2^n}+1\}$ is a couple if $2^{2^n}+1$ is a Fermat prime.
If $p$ is a prime and $p-1$ and $p+1$ are both practical, then we call $\{p-1,p,p+1\}$
a {\it sandwich of the first kind}. If $\{p,p+2\}$ is a twin prime pair and $p+1$ is practical, then we call $\{p,p+1,p+2\}$
a {\it sandwich of the second kind}. For example, $\{88,\, 89,\, 90\}$ is a sandwich of the first kind, while
$\{59,\,60,\,61\}$ is a sandwich of the second kind. See \cite[A210479]{S} for the list of the first 10,000 sandwiches of the first kind,
and \cite[A258838]{S} for the list of the first 10,000 sandwiches of the second kind.

\begin{conjecture}\label{Conj3.39} {\rm (2013-01-12) (i)} For any integer $n>8$ the interval $[n,2n]$ contains a sandwich of the first kind.

{\rm (ii)} For each $n=7,8,\ldots$ the interval $[n,2n]$ contains a sandwich of the second kind.

{\rm (iii)} For any integer $n>231$ the interval $[n,2n]$ contains four consecutive integers $p-1,p,p+1,p+2$ with $\{p,p+2\}$ a twin prime pair
and $\{p-1,p+1\}$ a twin practical pair.

{\rm (iv)} There are infinitely many quintuples $\{m-2,m-1,m,m+1,m+2\}$ with $\{m-1,m+1\}$ a twin prime pair and $m,m\pm2$ all practical.
\end{conjecture}
\begin{remark}\label{Rem3.39}. For those middle terms $m$ described in part (iv) of Conjecture \ref{Conj3.39}, the reader may consult \cite[A209236]{S}.
It is known that (cf. \cite{Me}) there are infinitely many practical numbers $m$ with $m\pm2$ also practical.
\end{remark}

\begin{conjecture}\label{Conj3.40} {\rm (i) (2013-01-23)} Each $n=4,5,\ldots$ can be written as $p+q$, where $\{p-1,p,p+1\}$ is a sandwich of the first kind, and $q$ is either prime or practical.

{\rm (ii) (2013-01-29)} Any integer $n>11$ can be written as $(1+\{n\}_2)p+q+r$, where $\{p-1,p,p+1\}$ and $\{q-1,q,q+1\}$
are sandwiches of the first kind, and $\{r-1,r,r+1\}$ is a sandwich of the second kind.
\end{conjecture}
\begin{remark}\label{Rem3.40} We have verified parts (i) and (ii) of Conjecture \ref{Conj3.40} for $n$ up to $10^8$ and $10^7$ respectively.
For numbers of representations related to parts (i) and (ii), see \cite[A210480 and A210681]{S}.
\end{remark}

\begin{conjecture}\label{Conj3.41} {\rm (2013-01-29) (i)} Any integer $n>6$ can be written as $p+q+r$ such that $\{p-1,p,p+1\}$ and $\{q-1,q,q+1\}$ are sandwiches of the first kind,
and $\{6r-1,6r,6r+1\}$ is a sandwich of the second kind.

{\rm (ii)} Every $n=3,4,\ldots$ can be expressed as $x+y+z$ with $x,y,z\in\Z^+$
such that $\{6x-1,6x,6x+1\}$, $\{6y-1,6y,6y+1\}$ and $\{6z-1,6z,6z+1\}$ are all sandwiches of the second kind.

{\rm (iiii)} Each integer $n>7$ can be written as $p+q+x^2$ with $x\in\Z$ such that $\{p-1,p,p+1\}$
is a sandwich of the first kind and $\{q-1,q,q+1\}$ is a sandwich of the second kind.
\end{conjecture}
\begin{remark}\label{Rem3.41} We also conjecture that each $n=3,4,\ldots$ can be written as the sum of two triangular numbers and a prime $p$ with $\{p-1,p,p+1\}$
a sandwich of the first kind. See \cite[A210681]{S} for related comments.
\end{remark}

\begin{conjecture}\label{Conj3.42} {\rm (2013-01-30) (i)} For any integer $n>8$, we can write $2n=p+2q+3r$, where $\{p-1,p,p+1\},\{q-1,q,q+1\}$ and $\{r-1,r,r+1\}$ are all sandwiches of the first kind.

{\rm (ii)} Each integer $n>5$ can be written as the sum of a prime $p$ with $\{p-1,p,p+1\}$ a sandwich of the first kind, a prime $q$ with $q+2$ also prime, and a Fibonacci number.
\end{conjecture}
\begin{remark}\label{Rem3.42} See \cite[A211190 and A211165]{S} for related data. We have verified part (ii) of Conjecture \ref{Conj3.42}
for $n$ up to $2,000,000$.
\end{remark}

\begin{conjecture}\label{Conj3.43} {\rm (i) (2013-01-14)} Any odd number $n>1$ can be expressed as $p+q$,
where $p$ is a Sophie Germain prime and $q$ is a practical number.

{\rm (ii) (2013-01-19)} For any integer $n>2$,
there is a practical number $q<n$ such that $n-q$ and $n+q$ are both prime or both practical.
\end{conjecture}
\begin{remark}\label{Rem3.43} We have verified this conjecture for $n$ up to $10^8$. See \cite[A209253 and A209312]{S} for related data.
We also conjecture that each positive integer can be represented as the sum of a practical number and a triangular number (cf. \cite[A208244]{S}),
which is an analog of the author's conjecture on sums of primes and triangular numbers (cf. \cite{S09}).
\end{remark}

\begin{conjecture}\label{Conj3.44} {\rm (2015-08-28) (i)} For any integer $n>6$, there is a prime $p<n$ such that
$n-(p+1)$ and $n+(p+1)$ are both prime or both practical.

{\rm (ii)} For any integer $n>2$ there is a prime $p<n$ such that $n-(p-1)$ and $n+(p-1)$ are both prime or both practical.
\end{conjecture}
\begin{remark}\label{Rem3.44} See \cite[A261653]{S} for related data, and compare this conjecture with Conjectures \ref{Conj2.2}, \ref{Conj2.3} and \ref{Conj3.43}.
\end{remark}

\begin{conjecture}\label{Conj3.45} {\rm (2015-07-12) (i)} There are infinitely many sandwiches $\{n-1,n,n+1\}$ of the first kind such that
$\{p_n-1,p_n,p_n+1\}$ is also a sandwich of the first kind.

{\rm (ii)} There are infinitely many sandwiches $\{n-1,n,n+1\}$ of the second kind such that
$\{p_n-1,p_n,p_n+1\}$ is a sandwich of the first kind.
\end{conjecture}
\begin{remark}\label{Rem3.45} See \cite[A257924 and A257922]{S} for related data.
\end{remark}

\section{On representations of positive rational numbers}
\setcounter{lemma}{0}
\setcounter{theorem}{0}
\setcounter{corollary}{0}
\setcounter{remark}{0}
\setcounter{equation}{0}
\setcounter{conjecture}{0}

It is well known that any positive rational number can be written as finitely many distinct unit fractions.
It is also known that the series $\sum_{n=1}^\infty1/p_n$ diverges as proved by Euler.

\begin{conjecture}\label{Conj4.1} {\rm (i) (2015-09-09)} For any positive rational number $r$, there are finitely many distinct primes
$q_1,\ldots,q_k$ such that
$$r=\sum_{j=1}^k\f1{q_j-1}.$$

{\rm (ii) (2015-09-12)} For any positive rational number $r$, there are finitely many distinct primes
$q_1,\ldots,q_k$ such that
$$r=\sum_{j=1}^k\f1{q_j+1}.$$

{\rm (iii) (2015-09-12)} For any positive rational number $r$, there are finitely many distinct practical numbers
$q_1,\ldots,q_k$ with $r=\sum_{j=1}^k1/q_j$.
\end{conjecture}
\begin{remark}\label{Rem4.1} For example,
$$2=\f1{2-1}+\f1{3-1}+\f1{5-1}+\f1{7-1}+\f1{13-1}=\f11+\f12+\f14+\f16+\f1{12}$$
with $2,3,5,7$ all prime and $1,2,4,6,12$ all practical, and
$$1=\f1{2+1}+\f1{3+1}+\f1{5+1}+\f1{7+1}+\f1{11+1}+\f1{23+1}$$
with $2,3,5,7,11,23$ all prime. Also,
\begin{align*}\f{10}{11}=&\f1{3-1}+\f1{5-1}+\f1{13-1}+\f1{19-1}+\f1{67-1}+\f1{199-1}
\\=&\f1{2+1}+\f1{3+1}+\f1{5+1}+\f1{7+1}+\f1{43+1}+\f1{131+1}+\f1{263+1}
\\=&\f12+\f14+\f1{8}+\f1{48}+\f1{132}+\f1{176}
\end{align*}
with $2,3,5,7,13,19,43,67,131,199,263$ all prime and $2,4,8,48,132,176$ all practical.
After learning Conjecture \ref{Conj4.1} from the author, Qing-Hu Hou verified parts (i) and (ii) in Nov. 2015 for all rational numbers $r\in(0,1)$ with denominators not exceeding 100.
The author would like to offer 500 US dollars as the prize for the first solution to parts (i) and (ii) of Conjecture \ref{Conj4.1}.
\end{remark}

\begin{conjecture}\label{Conj4.2} {\rm (2015-09-09)} Let $m$ be any positive integer.

{\rm (i)} All the rational numbers
     $$\sum_{i=j}^k\f1{(p_i-1)^m}\ \ \mbox{with}\  1\ls j \ls  k$$
are pairwise distinct! If
$$\sum_{i=j}^k\f1{(p_i-1)^m}\ \ \mbox{and}\ \ \sum_{r=s}^t \f1{(p_r-1)^m}$$
have the same fractional part with
$$0 < \min\{2,k\} \ls  j\ls  k,\  0 < \min\{2,t\} \ls s \ls t\ \mbox{and}\ j \ls s,$$
but the ordered pairs $(j,k)$ and $(s,t)$ are different, then we must have
$m = 1$ and
 $$\sum_{i=j}^k\f1{p_i-1} = 1 + \sum_{r=s}^t\f1{p_r-1};$$
moreover, either $(j,k) = (2,6)$ and $(s,t) = (5,5)$, or $(j,k) = (2,5)$ and $(s,t) = (18,18)$, or $(j,k) = (2,17)$ and $(s,t) =(6,18)$.

{\rm (ii)} If
$$\sum_{i=j}^k\f1{(p_i+1)^m}\ \ \mbox{and}\ \ \sum_{r=s}^t \f1{(p_r+1)^m}$$
have the same fractional part with
$$1\ls  j\ls  k,\  1\ls s \ls t\ \mbox{and}\ j \ls s,$$
but the ordered pairs $(j,k)$ and $(s,t)$ are different, then we must have
$m = 1$ and
$$\sum_{i=j}^k\f1{p_i+1} -\sum_{r=s}^t\f1{p_r+1}\in\{0,1\};$$
moreover, $(j,k) = (1,9)$ and $(s,t) = (6,8)$, or $(j,k) = (4,4)$ and $(s,t) = (8,10)$, or $(j,k) = (4,7)$ and $(s,t) =(5,10)$, or $(j,k) = (1,10)$ and $(s,t) = (5,7)$.

{\rm (iii)} For any integer $d > 1$, the rational numbers
$$\sum_{i=j}^k\f1{(p_i+d)^m}\ \ \mbox{with}\  1 \ls j \ls  k$$
have pairwise distinct fractional parts.
\end{conjecture}
\begin{remark}\label{Rem4.2} Recall that $\sum_{j=1}^\infty1/p_j$ diverges.
\end{remark}

Actually Conjecture \ref{Conj4.2} was motivated by our following conjecture whose proofs might involve primes.

\begin{conjecture}\label{Conj4.3} {\rm (i) (2015-09-09)} If $1/j+...+1/k$ and $1/s+...+1/t$ have the same fractional part with
$$0 < \min\{2,k\} \ls j\ls k,\ 0 < \min\{2,t\} \ls s\ls t\ \mbox{and}\  j\ls s,$$
but the ordered pairs $(j,k)$ and $(s,t)$ are different, then we have
      $$\f1j+\ldots+\f1k = 1+\f1s+\ldots+\f1t;$$
moreover, one of the following {\rm (a)-(d)} holds.

   {\rm (a)} $(j,k) = (2,6)$ and $(s,t) = (4,5)$,

   {\rm (b)} $(j,k) = (2,4)$ and $(s,t) = (12,12)$,

   {\rm (c)} $(j,k) = (2,11)$ and $(s,t) =(5,12)$,

   {\rm (d)} $(j,k) = (3,20)$ and $(s,t) = (7,19)$.

{\rm (ii) (2015-09-11)} Let $ a > b \gs 0$ and $m > 0$ be integers with $\gcd(a,b) = 1<\max\{a,m\}$.
Then the numbers
$$\sum_{i=j}^k\f1{(ai-b)^m}\ \ \mbox{with}\ 1 \ls j\ls k\ \mbox{and}\  (j > 1\ \mbox{if}\ k > a-b = 1)$$
have pairwise distinct fractional parts. Also, for each $r=0,1$, the numbers
$$\sum_{i=j}^k\f{(-1)^{i-jr}}{(ai-b)^m}\ \ \mbox{with}\ 1 \ls j\ls k\ \mbox{and}\  (j > 1\ \mbox{if}\ k > a-b = 1)$$
have pairwise distinct fractional parts.
\end{conjecture}
\begin{remark}\label{Rem4.3} In 1918 J. K\"urschak proved that for any integers $k \gs j>1$ the number $1/j+...+1/k$ is not an integer.
In 1946 P. Erd\H os and I. Niven \cite{EN} used Sylvester's theorem (which states that the product of $n$ consecutive integers greater than $n$ is divisible by a prime greater than $n$) to show that all the numbers $1/j+...+1/k$ with $1\ls j\ls k$ are pairwise distinct.
\end{remark}

If $d\in\Z^+$ is not a square, then the Pell equation $x^2-dy^2=1$ has infinitely many integral solutions.
Thus, for $r=a/b$ with $a,b\in\Z^+$ and $\gcd(a,b)=1$, if $r$ is not a square of rational numbers then there is a positive integer $k$ such that
$(ka)(kb)+1$ is a square, i.e., we can write $a/b=m/n$ with $m,n\in\Z^+$ such that $mn+1$ is a square.
Motivated by this, below we consider various representations of positive rational numbers.

\begin{conjecture}\label{Conj4.4} {\rm (i) (2015-07-03)} The set
$$\l\{\f mn:\ m,n\in\Z^+\ \mbox{and}\ p_m+p_n\ \t{is a square}\r\}$$
contains any positive rational number $r$. Also, any rational number $r>1$ can be written as $m/n$ with $m,n\in\Z^+$ such that $p_m-p_n$ is a square.

{\rm (ii) (2015-08-20)} Any positive rational number $r\not=1$ can be written as $m/n$ with $m,n\in\Z^+$  such that $p_{p_m}+p_{p_n}$ is a square.
\end{conjecture}
\begin{remark}\label{Rem4.4} We have verified part (i) of Conjecture \ref{Conj4.4} for all those rational numbers $r=a/b$ with $a,b\in\{1,\ldots,200\}$
(cf. \cite[A259712 and A257856]{S}) and part (ii) of Conjecture \ref{Conj4.4} for all those rational numbers $r=a/b\not=1$ with $a,b\in\{1,\ldots,60\}$.
For example, $2=20/10$ with $p_{20}+p_{10}=71+29=10^2$, and $2=92/46$ with
$p_{p_{92}}+p_{p_{46}}=p_{479}+p_{199}=3407+1217=68^2$.
\end{remark}

\begin{conjecture}\label{Conj4.5} {\rm (2015-07-08)} The set
$$\l\{\f mn:\ m,n\in\Z^+,\ \t{and}\ \varphi(m)\ \t{and}\ \sigma(n)\ \t{are both squares}\r\}$$
contains any positive rational number $r$.
\end{conjecture}
\begin{remark}\label{Rem4.5} We have verified Conjecture \ref{Conj4.5} for all those $r=a/b$ with $a,b\in\{1,\ldots,150\}$ (cf. \cite[A259915 and A259916]{S}).
For example, $4/5 = 136/170$ with $\varphi(136) = 8^2$ and $\sigma(170) = 18^2$, and $5/4 = 1365/1092$ with $\varphi(1365) = 24^2$ and $\sigma(1092) = 56^2$.
\end{remark}

\begin{conjecture}\label{Conj4.6} {\rm (i) (2015-07-05)} Any positive rational number $r$ can be written as
$m/n$ with $m,n\in\Z^+$ such that $\pi(m)\pi(n)$ is a positive square.

{\rm (ii) (2015-07-06)} Any positive rational number $r$ can be written as
$m/n$ with $m,n\in\Z^+$ such that $\pi(m)$ and $\pi(\pi(n))$ are positive squares.
\end{conjecture}
\begin{remark}\label{Rem4.6} We have verified part (i) of this conjecture for all those rational numbers $r=a/b$ with $a,b\in\{1,\ldots,60\}$. See \cite[A259789]{S} for related data.
For example, $49/58 = 1076068567/1273713814$ with
$$\pi(1076068567)\pi(1273713814) = 54511776\cdot63975626 = 59054424^2.$$
\end{remark}

\begin{conjecture}\label{Conj4.7} {\rm (2015-07-10)} Each positive rational number $r<1$ can be written as $m/n$ with $1<m<n$ such that $\pi(m)^2+\pi(n)^2$ is a square.
Also, any rational number $r>1$ can be written as $m/n$ with $m>n>1$ such that $\pi(m)^2-\pi(n)^2$ is a square.
\end{conjecture}
\begin{remark}\label{Rem4.7} We have verified this conjecture for all those rational numbers $r=a/b$ with $a,b\in\{1,\ldots,50\}$. See \cite[A255677]{S} for related data.
For example, $23/24 = 19947716/20815008$ with
$$\pi(19947716)^2 + \pi(20815008)^2 = 1267497^2 + 1319004^2 = 1829295^2,$$
and $7/3 = 26964/11556$ with
$$\pi(26964)^2 - \pi(11556)^2 = 2958^2 - 1392^2 = 2610^2.$$
\end{remark}

Motivated by Conjecture \ref{Conj4.7}, we raise the following conjecture which sounds interesting and challenging.

\begin{conjecture}\label{Conj4.8} {\rm (i) (2015-07-11)} For any $n\in\Z^+$, there are distinct primes $p,q,r$ such that
$\pi(pn)^2=\pi(qn)^2+\pi(rn)^2$.

{\rm (ii) (2015-07-13)} For any $n\in\Z^+$, there are distinct primes $p,q,r$ with $\pi(pn)=\pi(qn)\pi(rn)\ (or\ \pi(pn)=\pi(qn)+\pi(rn))$.
\end{conjecture}
\begin{remark}\label{Rem4.8} See \cite[A257364 and A257928]{S} for related data.
\end{remark}

\begin{conjecture}\label{Conj4.9} {\rm (i) (2015-07-02)} Any positive rational number $r$ can be written as $m/n$ with $m,n\in\Z^+$ such that $p(m)^2+p(n)^2$ is prime, where
$p(\cdot)$ is the partition function.

{\rm (ii) (2015-08-20)} Any positive rational number $r\not=1$ can be written as $m/n$ with $m,n\in\Z^+$ such that $p(p_m)+p(p_n)$ is prime.
\end{conjecture}
\begin{remark}\label{Rem4.9}  Conjecture \ref{Conj4.9} implies that there are infinitely many primes of the form $p(m)^2+p(n)^2$ with $m,n\in\Z^+$ as well as primes of the form $p(q)+p(r)$ with $q$ and $r$ both prime. We have verified part (i) for all those rational numbers $r=a/b$ with $a,b\in\{1,\ldots,100\}$,
and part (ii) for all those rational numbers $r=a/b\not=1$ with $a,b\in\{1,\ldots,37\}$. See \cite[A259531, A259678, A261513 and A261515]{S} for related data.
For example, $4/5 = 124/155$ with
\begin{align*} p(124)^2 + p(155)^2 =& 2841940500^2 + 66493182097^2
\\=& 4429419891190341567409
\end{align*}
prime, and $3=138/46$ with
\begin{align*}p(p_{138}) + p(p_{46}) = &p(787)+p(199)
\\=&3223934948277725160271634798+3646072432125
\\=& 3223934948277728806344066923
\end{align*}
prime.
\end{remark}

\begin{conjecture}\label{Conj4.10} {\rm (2015-08-17)} Any positive rational number $r$ can be written as $m/n$, where $m$ and $n$ are positive integers with
$(m\pm1)^2+n^2$ and $m^2+(n\pm1)^2$ all prime.
\end{conjecture}
\begin{remark}\label{Rem4.10} We have verified this for all those $r=a/b$ with $a,b\in\{1,\ldots,60\}$. See \cite[A261382]{S} for related data.
It is easy to prove that if $m$ and $n$ are positive integers with $(m\pm1)^2+n^2$ and $m^2+(n\pm1)^2$ all prime then either $m=n=2$
or $m\eq n\eq0\pmod5$.
\end{remark}

\begin{conjecture}\label{Conj4.11} {\rm (i) (2015-06-28)} Each rational number $r>0$ can be written as $m/n$, where $m$ and $n$ are positive integers with
$$p_m\pm m,\ p_n\pm n,\ p_m+n\ \mbox{and}\ p_n+m$$
all prime.

{\rm (ii) (2015-07-02)} Any rational number $r>0$ can be written as $m/n$, where $m$ and $n$ are positive integers with
$$m^2+p_m^2,\ n^2+p_n^2,\ m^2+p_n^2\ \mbox{and}\ n^2+p_m^2$$
all prime.

{\rm (iii) (2015-08-15)} Any rational number $r>0$ can be written as $m/n$ with $m$ and $n$ in the set
\begin{align*}&\{k\in\Z^+:\ k+1,\ k^2+1\ \mbox{and}\ k^2+p_k^2\ \mbox{are all prime}\}
\\=&\{q-1:\ q,\ (q-1)^2+1\ \mbox{and}\ (q-1)^2+p_{q-1}^2\ \mbox{are all prime}\}.
\end{align*}
\end{conjecture}
\begin{remark}\label{Rem4.11} We have verified parts (i)-(ii) for those $r=a/b$ with $a,b\in\{1,\ldots,150\}$
and part (iii) for those $r=a/b$ with $a,b\in\{1,\ldots,60\}$. See \cite[A259492 and A261339]{S} for related data.
\end{remark}

\begin{conjecture}\label{Conj4.12} {\rm (i) (2015-06-30)} Let
$$U:=\{n\in\Z^+:\ n\pm1\ \t{and}\ p_n+2\ \t{are all prime}\}.$$
Then any positive rational number $r$ can be written as $m/n$ with $m,n\in U$.

{\rm (ii) (2015-06-28)} Let
$$V:=\{n\in\Z^+:\ p_n+2\ \t{and}\ p_{p_n}+2\ \t{are both prime}\}.$$
Then any positive rational number $r$ can be written as $m/n$ with $m,n\in V$.

{\rm (iii) (2015-06-12)} Let
$$Q:=\{q\in\Z^+:\ q\ \t{is practical with}\ q\pm1\ \t{twin prime}\}.$$
Then any positive rational number $r$ can be written as $q/q'$ with $q,q'\in Q$.
\end{conjecture}
\begin{remark}\label{Rem4.12} We have verified part (i) for all those $r=a/b$ with $a,b\in\{1,\ldots,100\}$,  part (ii) for all those
$r=a/b$ with $a,b\in\{1,\ldots,400\}$, and part (iii) for all those
$r=a/b$ with $a,b\in\{1,\ldots,1000\}$. See \cite[A259539, A259540, A259487, A259488 and A258836]{S} for related data.
For example, $4/5 = 11673840/14592300$ with 11673840 and 14592300 in the set $U$.
\end{remark}

Motivated by part (i) of Conjecture \ref{Conj4.12} and \cite[Conjecture 3.7(i)]{S15}, we pose the following conjecture.
\begin{conjecture}\label{Conj4.13} {\rm (2015-07-01)} There are infinitely many positive integers $n$ such that the seven numbers
$$n\pm1,\ p_n+2,\ p_n\pm n,\ np_n\pm1$$
are all prime.
\end{conjecture}
\begin{remark}\label{Rem4.13} We have listed the first 160 such positive integers $n$ the least of which is 2523708 (cf. \cite[A259628]{S}).
\end{remark}

\begin{conjecture}\label{Conj4.14} {\rm (2015-08-24)} Any positive rational number $r$ can be written as $m/n$, where $m$ and $n$ belong to the set
$$\{k\in\Z^+:\ p_k+2,\ p_k+6\ \mbox{and}\ p_k+8\ \mbox{are all prime}\}.$$
Also, each positive rational number $r$ can be written as $m/n$, where $m$ and $n$ belong to the set
$$\{k\in\Z^+:\ p_k+4,\ p_k+6\ \mbox{and}\ p_k+10\ \mbox{are all prime}\}.$$
\end{conjecture}
\begin{remark}\label{Rem4.14}  This conjecture implies that there are infinitely many prime quadruples $(p,p+2,p+6,p+8)$ as well as
$(p,p+4,p+6,p+10)$, which is a special case of Schinzel's Hypothesis. See \cite[A261541]{S} for related data. For example,
$3/4=m/n$ with $m=20723892$ and $n=27631856$, and
\begin{align*}p_{m}+2&=387875563,\ p_{m}+6=387875567,\ p_{m}+8=387875569,
\\p_{n}+2&=525608593,\ p_{n}+6=525608597,\ p_{n}+8=525608599
\end{align*}
are all prime.
\end{remark}

\begin{conjecture}\label{Conj4.15} {\rm (2015-08-23)} Any positive rational number $r$ can be written as $m/n$, where $m$ and $n$ belong to the set
$$W=\{k\in\Z^+:\ p_k+2\ \mbox{is prime}\ \mbox{and}\ p_{p_k+2}-p_{p_k}=6\}.$$
\end{conjecture}
\begin{remark}\label{Rem4.15}  See \cite[A261528 and A261533]{S} for related data. For example, $2=1782/891$ with $891$ and $1782$ in the set $W$.
Conjecture \ref{Conj4.15} implies that there are infinitely many twin prime pairs $\{q,q+2\}$ with $p_{q+2}-p_q=6$.
\end{remark}

\begin{conjecture}\label{Conj4.16} {\rm (2015-08-14)} Each positive rational number $r$ can be written as $m/n$ with $m$ and $n$ in the set
$$\{k\in\Z^+:\ p_k^2-2\ \mbox{and}\ p_{p_k}^2-2\ \mbox{are both prime}\}.$$
\end{conjecture}
\begin{remark}\label{Rem4.16} We have verified this for all those $r=a/b$ with $a,b\in\{1,\ldots,300\}$.
See \cite[A261281]{S} for related data.
\end{remark}

\begin{conjecture}\label{Conj4.17} {\rm (i) (2014-05-14)} For any prime $p>5$, there is a positive square $k^2<p$ such that the inverse of $k^2$ modulo $p$ is prime,
where the inverse of $a\in\{1,\ldots,p-1\}$ modulo $p$ denotes the unique $x\in\{1,\ldots,p-1\}$ with $ax\eq1\pmod{p}$.

{\rm (ii) (2015-08-18)} Any positive rational number $r\ls1$ can be written as $m/n$ with $m,n\in\Z^+$ such that the inverse of $m$ modulo $p_n$ is a square.
\end{conjecture}
\begin{remark}\label{Rem4.17} We have checked part (i) of Conjecture \ref{Conj4.17} for those primes $p<1.8\times10^8$.
See \cite[A242425 and A242441]{S} for related data. For example, the inverse of $4^2$ modulo 23 is the prime 13.
\end{remark}

\begin{conjecture}\label{Conj4.18} {\rm (2014-08-26) (i)} Any integer $n>2$ with $n\not=8$ can be written as $k+m$ with $k,m\in\Z^+$ and $k\not=m$
such that $p_k$ is a primitive root modulo $p_m$ and $p_m$ is also a primitive root modulo $p_k$.

{\rm (ii)} Any positive rational number $r\not=1$ can be written as $m/n$ with $m,n\in\Z^+$ such that $p_m$ is a primitive root modulo $p_n$
and also $p_n$ is a primitive root modulo $p_m$.
\end{conjecture}
\begin{remark}\label{Rem4.18} See \cite[A261387]{S} for related data and comments.
\end{remark}

\begin{conjecture}\label{Conj4.19} {\rm (2015-07-20)} Let $n\in\Z^+$ and $s,t\in\{1,-1\}$. Then any positive rational number $r_0$ can be written as
$(p_{qn}+s)/(p_{rn}+t)$ with $q$ and $r$ both prime, unless $n>r_0=1$ and $\{s,t\}=\{1,-1\}$.
\end{conjecture}
\begin{remark}\label{Rem4.19}  We have verified this conjecture in the case $n=1$ for all those $r_0=a/b$ with $a,b\in\{1,\ldots,500\}$ (cf. \cite[A258803]{S}).
For $n=2,\ldots,10$ we have verified Conjecture \ref{Conj4.19} for all those $r_0=a/b$ with $a,b\in\{1,\ldots,30\}$ (cf. \cite[A260252]{S}). For example,
$23 = (p_{17209} - 1)/(p_{1039} - 1) = (190579 - 1)/(8287 - 1)$ with 1039 and 17209 both prime.
\end{remark}

\begin{conjecture}\label{Conj4.20} {\rm (2015-08-02) (i)} If $a,b,c$ are positive integers with $\gcd(a,b) = \gcd(a,c) = \gcd(b,c) = 1$, and $a\not=b$ and $a+b\eq c\pmod 2$,
then for any $n\in\Z^+$ the linear equation $ax-by = c$ has solutions with $x$ and $y$ in the set $\{p_{qn}: q\ \mbox{is prime}\}$.

{\rm (ii)} Let $a$ and $b$ be relatively prime positive integers, and let $c$ be any integer. For any $n\in\Z^+$, the linear equation $ax-by = c$ has solutions with $x$ and $y$ in the set
 $\{\pi(pn):\ p\ \mbox{is prime}\}$.
\end{conjecture}
\begin{remark}\label{Rem4.20} Note that part (i) of Conjecture \ref{Conj4.20} is an extension of Conjecture \ref{Conj4.19}. In the
$a = c = 1$ and $b = 2$, it asserts that for any $n\in\Z^+$ there are primes $q$ and $r$ such that
$2p_{qn}+1 = p_{rn}$. This implies that there are infinitely many Sophie Germain primes.
Also, part (ii) of Conjecture \ref{Conj4.20}
with $c = 0$ asserts that for any $n\in\Z^+$ the set
$$\left\{\f{\pi(pn)}{\pi(qn)}:\ p\ \mbox{and}\ q\ \mbox{are primes}\right\}$$ contains all positive rational numbers
(cf. \cite[A260232]{S}). We have checked both parts of the conjecture for $a,b,c = 1,\ldots,20$ and  $n = 1,\ldots,30$. For related data, see
\cite[A260886 and A260888]{S}.
\end{remark}

Recall that a prime $p$ is called a Chen prime if $p+2$ is a product of at most two primes.
In 1973 J. Chen \cite{C} proved that there are infinitely many Chen primes.

\begin{conjecture}\label{Conj4.21} {\rm (i) (2015-07-14)} For any positive integer $n$, there are $i,j,k\in\Z^+$ with $i\not=j$ such that $p_{kn}+2=p_{in}p_{jn}$.

{\rm (ii) (2015-07-15)} For any positive integer $n$, there are $i,j,k\in\Z^+$ with $i\not=j$ such that $p_{kn}^2-2=p_{in}p_{jn}$.
\end{conjecture}
\begin{remark}\label{Rem4.21} See \cite[A257926 and A260080]{S} for related data.
Clearly, part (i) of Conjecture \ref{Conj4.21} implies that there are infinitely many Chen primes.
\end{remark}

\begin{conjecture}\label{Conj4.22} {\rm (2015-07-15)} Let $d$ be a nonzero integer and let $n\in\Z^+$. Set
$$D:=\{p_{kn}+d:\ k=1,2,3,\ldots\}.$$

{\rm (i)} If $\gcd(6,d)=1$, then there are two distinct elements $x$ and $y$ of $D$ with $x+y\in D$ and $x-y\in D$.

{\rm (ii)} For each $k=1,2$, we have $xy=z^k$ for some distinct elements $x,y,z$ of $D$.
\end{conjecture}
\begin{remark}\label{Rem4.22} See \cite[A260078, A257938 and A260082]{S} for related data.
\end{remark}

\begin{conjecture}\label{Conj4.23} {\rm (2015-07-17) (i)} Let $a,n\in\Z^+$ and $b,c\in\Z$ with $\gcd(a,b,c)=1$, $2\nmid(a+b+c)$ and $3\nmid\gcd(b,a+c)$.
If $b^2-4ac$ is not a square, then there are $x,y\in\{p_{kn}:\ k=1,2,3,\ldots\}$
such that $y=ax^2+bx+c$.

{\rm (ii)} For any $a,n\in\Z^+$ and $b,c\in\Z$, there are $x,y\in\{\pi(pn):\ p \ \mbox{is prime}\}$
such that $y=ax^2+bx+c$.
\end{conjecture}
\begin{remark}\label{Rem4.23} See \cite[A260120 and A260140]{S} for related data. Part (i) of Conjecture \ref{Conj4.23} implies that
for any $n\in\Z^+$ there are $j,k\in\Z^+$ with $p_{kn}^2-2=p_{jn}$ (or $(p_{kn}-1)^2=p_{jn}-1$).
Part (ii) of Conjecture \ref{Conj4.23} implies that for any $n\in\Z^+$ there are primes $p$ and $q$ with $\pi(pn)=\pi(qn)^2.$
\end{remark}

\begin{conjecture}\label{Conj4.24} {\rm (2015-08-14)} Let
$$S_1:=\{q+1:\ q\ \mbox{and}\ p_q+2\ \mbox{are both prime}\}$$
and
$$S_2:=\{q-1:\ q\ \mbox{and}\ p_q-2\ \mbox{are both prime}\}.$$
For any $i,j\in\{1,2\}$, each positive rational number $r$ can be written as $m/n$ with $m\in S_i$ and $n\in S_j$, unless $i\not=j$ and $r=1$.
\end{conjecture}
\begin{remark}\label{Rem4.24} See \cite[A261295]{S} for related data. For example,
$4/5=15648/19560$ with $15647,\ p_{15647}+2=171763,\ 19559$ and $p_{19559}+2=219409$ all prime.
A twin prime pair $\{p,p+2\}$ with $\pi(p)$ also prime is called a super twin prime pair (cf. \cite [Conjecture 3.2 and Remark 3.2]{S15}).
\end{remark}

\begin{conjecture}\label{Conj4.25} {\rm (2015-08-18)} Let $s,t\in\{1,-1\}$. Then any positive rational number $r$ can be written as $m/n$ with $m$ and $n$ in the set
$$K_{s,t}:=\{k\in\Z^+:\ p_{p_k}+sp_k+t=p_q\ \ \mbox{for some prime}\ q\}.$$
\end{conjecture}
\begin{remark}\label{Rem4.25} This implies that for any $s,t\in\{\pm1\}$ there are infinitely many primes $q$ with $p=p_q+sq+t$ and $\pi(p)$ both prime.
See \cite[A260753 and A261136]{S} for related data. For example, $3=6837/2279$, and
$$p_{p_{6837}}-p_{6837}+1 = p_{68777}-68777+1 = 865757-68776 = 796981 = p_{63737}$$
with 63737 prime, and
$$p_{p_{2279}}-p_{2279}+1 = p_{20147}-20147+1 = 226553-20146 = 206407 = p_{18503}$$
with $18503$ prime.
\end{remark}

\begin{conjecture}\label{Conj4.26} {\rm (2015-08-16)} Any positive rational number can be written as $m/n$, where $m$ and $n$
are positive integers with $p_{p_m}p_{p_n}=p_q+2$ for some prime $q$.
\end{conjecture}
\begin{remark}\label{Rem4.26}  See \cite[A261352 and A261353]{S} for related data. For example, $4=2424/606$ and
$$p_{p_{2424}}p_{p_{606}} = p_{21589}p_{4457} = 244471\cdot42643 = 10424976853 = p_{473490161}+2$$
with 473490161 prime. Conjecture \ref{Conj4.26} implies that there are infinitely many prime triples $(q,r,s)$ with $p_q+2=p_rp_s$.
\end{remark}

\begin{conjecture}\label{Conj4.27} {\rm (2014-08-17) (i)} Let $d$ be any nonzero integer. Then
any positive rational number $r$ can be written as $m/n$ with $m,n\in\Z^+$ such that
$(p_{p_m}+d)(p_{p_n}+d)=p_{q}+d$ for some prime $q$.

{\rm (ii)} For any nonzero integer $d$, there are infinitely many prime triples $(q,r,s)$ with $q,r,s$ distinct such that
$(p_q+d)^2=(p_r+d)(p_s+d)$.
\end{conjecture}
\begin{remark}\label{Rem4.27} See \cite[A261385 and A261395]{S} for related data and comments. Clearly, for each $d\in\Z\sm\{0\}$,  part (i) of Conjecture
\ref{Conj4.27} implies that the equation $xy=z$ has infinitely many solutions with $x,y,z\in\{p_q+d:\ q\ \mbox{is prime}\}$,
and part (ii) of Conjecture \ref{Conj4.27} implies that the set $\{p_q+d:\ q\ \mbox{is prime}\}$ contains infinitely many nontrivial 3-term geometric progressions.
\end{remark}

\begin{conjecture}\label{Conj4.28} {\rm (2015-08-16) (i)} Let $a,b,c\in\Z^+$ with $a\not=b$, $a+b\eq c\pmod 2$ and $\gcd(a,b)=\gcd(a,c)=\gcd(b,c)=1$.
Then any positive rational number $r$ can be written as $m/n$ with $m$ and $n$ in the set
$$\{k\in\Z^+:\ ap_q-bp_{p_k}=c\ \mbox{for some prime}\ q\},$$
thus there are infinitely many pairs of primes $q$ and $r$ such that $ap_q-bp_r=c$.

{\rm (ii)} Let $a\in\Z^+$ and $b,c\in\Z$ with $\gcd(a,b,c)=1$. If $2\nmid (a+b+c)$, $3\nmid\gcd(b,a+c)$, and $b^2-4ac$ is not a square, then the equation $y=ax^2+bx+c$
has infinitely many solutions with $x,y\in\{p_q:\ q\ \mbox{is prime}\}$.
\end{conjecture}
\begin{remark}\label{Rem4.28} See \cite[A261361, A261362 and A261354]{S} for related data and comments.
Clearly, part (i) of Conjecture \ref{Conj4.28} implies that there are infinitely many prime pairs $q$ and $r$ with $2p_q+1=p_r$, and
part (ii) of Conjecture \ref{Conj4.28} implies that there are infinitely many prime pairs $q$ and $r$ with $p_q^2-2=p_r$.
\end{remark}

\begin{conjecture}\label{Conj4.29} {\rm (i) (2015-08-18)}
For any $j=\pm1$ and $n\in\Z^+$, there is a positive integer $k$ such that $kn+j = p_q$ and $k^2n+1 =p_r$ for some pair of primes $q$ and $r$.

{\rm (ii) (2015-08-20)} Each positive rational number $r\ls 1$ can be written as $m/n$, where $m$ and $n$ are positive integers such that
$p_{p_m},p_{p_n},p_{p_k},p_{p_l}$ form a four-term arithmetic progression for some $k,l\in\Z^+$.

{\rm (iii) (2015-08-25)} Any positive rational number $r$ can be written as $m/n$, where $m$ and $n$ are positive integers
with $(p_{p_{p_m}}+p_{p_{p_n}})/2=p_{p_q}$ for some prime $q$.
\end{conjecture}
\begin{remark}\label{Rem4.29} See \cite[A261437, A261462 and A261583]{S} for related data.
\end{remark}

Motivated by Conjecture \ref{Conj4.29}, we define
$$p^{(1)}_n=p_n,\ \mbox{and}\ p^{(m+1)}_n=p^{(m)}_{p_n}\quad\mbox{for}\ m,n=1,2,3,\ldots,$$
and pose the following conjecture.

\begin{conjecture}\label{Conj4.30} {\rm (2015-08-25) (i)} If $q\in\Z^+$ and $a\in\Z$ are relatively prime, then for any $m\in\Z^+$
there are infinitely many $n\in\Z^+$ with $p^{(m)}_n\eq a\pmod{q}$.

{\rm (ii)} For any integer $k>2$ and $m>0$, the set $P_m:=\{p^{(m)}_n:\ n\in\Z^+\}$ contains infinitely many nontrivial $k$-term arithmetic progressions.

{\rm (iii)} For any $m,n\in\Z^+$, we have
$$\f{\root{n+1}\of{p^{(m+1)}_{n+1}}}{\root{n}\of{p^{(m+1)}_n}}
<\f{\root{n+1}\of{p^{(m)}_{n+1}}}{\root{n}\of{p^{(m)}_n}}<1.$$
\end{conjecture}
\begin{remark}\label{Rem4.30} Part (i) of Conjecture \ref{Conj4.30} is an extension of Dirichlet's theorem on primes in arithmetic progressions, and part (ii)
of Conjecture \ref{Conj4.30} is an extension of the Green-Tao theorem \cite{GT}. Part (iii) is an analog of Firoobakht's conjecture (cf. \cite{S13a}),
we also conjecture that the sequence $(\root{n}\of{q_n})_{n\gs 3}$ is strictly decreasing if $q_n$ denotes the $n$-th practical number.
\end{remark}

\medskip

\Ack.  The work was supported by the National Natural Science Foundation of China (Grant No. 11571162).

\end{document}